\documentclass[11pt]{article}

\usepackage{tikz}
\usepackage{subfigure}
\usepackage{scalefnt}

\usepackage{amsmath,amsfonts,amsthm,enumerate}

\newtheorem{defn}{{\bf Definition}}[section]
\newtheorem{eg}[defn]{{\bf Example}}
\newtheorem{lemma}[defn]{{\bf Lemma}}

\newtheorem{theo}[defn]{{\bf Theorem}}
\newtheorem{cor}[defn]{{\bf Corollary}}
\newtheorem{remark}[defn]{{\bf Remark}}
\newtheorem{conj}[defn]{{\bf Conjecture}}
\newtheorem{qn}[defn]{{\bf Question}}

\begin{document}

\title{\bf On \boldmath{$k$}-stellated and \boldmath{$k$}-stacked spheres }
\author{{\bf Bhaskar Bagchi}$^{\rm a}$ and {\bf Basudeb
Datta}$^{\rm b}$ }

\date{}

\maketitle

\vspace{-5mm}

\noindent {\small $^{\rm a}$Theoretical Statistics and Mathematics Unit, Indian Statistical Institute, Bangalore
560\,059, India.
\newline  \mbox{} \hspace{.2mm} bbagchi@isibang.ac.in

\smallskip

\noindent $^{\rm b}$Department of Mathematics, Indian Institute of Science, Bangalore 560\,012,  India. \newline
\mbox{} \hspace{.2mm}   dattab@math.iisc.ernet.in}

\date{}

\maketitle

\begin{center}

\date{May 16, 2013}

\end{center}

\vspace{-2mm}

\hrule

\medskip

 \centerline{\sc Abstract}

\medskip

{\small We introduce the class $\Sigma_k(d)$ of $k$-stellated (combinatorial) spheres of dimension $d$ ($0 \leq
k\leq d+1$) and compare and contrast it with the class ${\mathcal S}_k(d)$ ($0\leq k\leq d$) of $k$-stacked
homology $d$-spheres. We have $\Sigma_1(d) = {\mathcal S}_1(d)$, and $\Sigma_k(d) \subseteq {\mathcal S}_k(d)$
for $d\geq 2k-1$. However, for each $k \geq 2$ there are $k$-stacked spheres which are not $k$-stellated. For
$d\leq 2k -2$, the existence of $k$-stellated spheres which are not $k$-stacked remains an open question.

We also consider the class ${\mathcal W}_k(d)$ (and ${\mathcal K}_k(d)$) of simplicial complexes all whose
vertex-links belong to $\Sigma_k(d-1)$ (respectively, ${\mathcal S}_k(d-1)$). Thus, ${\mathcal W}_k(d) \subseteq
{\mathcal K}_k(d)$ for $d\geq 2k$, while ${\mathcal W}_1(d) = {\mathcal K}_1(d)$. Let $\overline{{\mathcal
K}}_k(d)$ denote the class of $d$-dimensional complexes all whose vertex-links are $k$-stacked balls. We show
that for $d\geq 2k+2$, there is a natural bijection $M\mapsto \overline{M}$ from ${\mathcal K}_k(d)$ onto
$\overline{{\mathcal K}}_k(d+1)$ which is the inverse to the boundary map $\partial \colon \overline{{\mathcal
K}}_k(d+1) \to {\mathcal K}_k(d)$.

Finally, we complement the tightness results of our recent paper \cite{bd16} by showing that, for any field
$\mathbb{F}$, an $\mathbb{F}$-orientable $(k+1)$-neighborly member of ${\mathcal W}_k(2k+1)$ is
$\mathbb{F}$-tight if and only if it is $k$-stacked. }

\bigskip

{\small


\smallskip

 }

\medskip

\hrule

\section{Introduction}

By a homology sphere/ball, we mean an $\mathbb{F}$-homology sphere/ball for some field $\mathbb{F}$. In this
paper, we introduce the class $\Sigma_k(d)$, $0\leq k \leq d +1$, of {$k$-stellated} triangulated $d$-spheres and
compare it with the class ${\mathcal S}_k(d)$, $0\leq k \leq d$, of {$k$-stacked} homology $d$-spheres. We have
the filtration
$$
\Sigma_0(d) \subseteq \Sigma_1(d) \subseteq \cdots \subseteq \Sigma_{d}(d) \subseteq \Sigma_{d +1}(d)
$$
of the class of all combinatorial $d$-spheres, and the comparable filtration
$$
{\mathcal S}_0(d) \subseteq {\mathcal S}_1(d) \subseteq \cdots \subseteq {\mathcal S}_{d}(d)
$$
of the class of all homology $d$-spheres. The {\em standard $d$-sphere} $S^{\,d}_{d + 2}$ is the unique $(d +
2)$-vertex triangulation of the $d$-sphere. It may be described as the boundary complex of the $(d +
1)$-dimensional geometric simplex. The standard sphere $S^{\,d}_{d + 2}$ is the unique member of $\Sigma_0(d) =
{\mathcal S}_0(d)$. We also have the equality $\Sigma_1(d) = {\mathcal S}_1(d)$. In the existing literature, the
members of ${\mathcal S}_1(d)$ are known as the $d$-dimensional {\em stacked spheres}. For $d\geq 2k - 1$, we
have the inclusion $\Sigma_k(d) \subseteq {\mathcal S}_k(d)$. However, for each $k \geq 2$, there are $k$-stacked
spheres which are not $k$-stellated.

In parallel with these classes of homology spheres, we also consider the classes $\widehat{\Sigma}_k(d)$ and
$\widehat{\mathcal S}_k(d)$ of {$k$-shelled} $d$-balls and {$k$-stacked} homology $d$-balls, respectively. We
have the filtration
$$
\widehat{\Sigma}_0(d) \subseteq \widehat{\Sigma}_1(d) \subseteq \cdots \subseteq \widehat{\Sigma}_{d}(d)
$$
of the class of all shellable $d$-balls, and the comparable filtration
$$
\widehat{\mathcal S}_0(d) \subseteq \widehat{\mathcal S}_1(d) \subseteq \cdots \subseteq \widehat{\mathcal
S}_{d}(d)
$$
of the class of all homology $d$-balls. The {\em standard $d$-ball} $B^{\,d}_{d + 1}$ is the unique $(d +
1)$-vertex triangulation of the $d$-dimensional ball. It may be described as the face complex of the
$d$-dimensional geometric simplex. The standard ball $B^{\,d}_{d + 1}$ is the unique member of
$\widehat{\Sigma}_0(d)= \widehat{\mathcal S}_0(d)$. We also have the equality $\widehat{\Sigma}_1(d) =
\widehat{\mathcal S}_1(d)$ and for all $d \geq k$ we have the inclusion $\widehat{\Sigma}_k(d) \subseteq
\widehat{\mathcal S}_k(d)$. However, for each $k \geq 2$, there are $k$-stacked balls which are not $k$-shelled.

While a {\em $k$-stellated $d$-sphere} is defined as a triangulated $d$-sphere which may be obtained from
$S^{\,d}_{d + 2}$ by a finite sequence of bistellar moves of index $< k$, a {\em $k$-shelled $d$-ball} is a
triangulated $d$-ball obtained from $B^{\,d}_{d+1}$ by a finite sequence of shelling moves of index $< k$. A {\em
$k$-stacked homology $d$-ball} is a homology $d$-ball all whose faces of codimension $k+1$ (i.e., dimension
$d-k-1$) are in its boundary. A {\em $k$-stacked homology $d$-sphere} is a homology $d$-sphere which may be
represented as the boundary of a $k$-stacked $(d+1)$-ball. The boundary of any $k$-shelled $(d+1)$-ball is a
$k$-stellated $d$-sphere. Conversely, when $d\geq 2k-1$, any $k$-stellated $d$-sphere may be represented as the
boundary of a $k$-shelled $(d+1)$-ball. A homology ball is $k$-shelled if and only if it is $k$-stacked and
shellable. Each $k$-stacked homology (respectively $k$-shelled) ball is the antistar of a vertex in a $k$-stacked
homology (respectively $k$-stellated) sphere. Murai and Nevo \cite{mn} proved that, when $d\geq 2k$, for any
$k$-stacked homology $d$-sphere $S$, there is a unique $k$-stacked homology $(d+1)$-ball $\overline{S}$ whose
boundary is $S$. The ball $\overline{S}$ has a natural and intrinsic description in terms of the combinatorics of
$S$. We point out the Murai-Nevo theorem is an immediate consequence of the following lemma\,: ``In any
$k$-stacked homology ball, all the missing faces have dimension $\leq k$". As another consequence of this lemma,
we show that the dimension $t$ of any missing face in a $k$-stacked homology $d$-sphere ($d\geq 2k+1$) satisfies
$t\leq k$ or $t \geq d-k+1$.

We consider the class ${\mathcal W}_k(d)$, $0\leq k\leq d$ (and ${\mathcal K}_k(d)$, $0\leq k\leq d-1$) of
simplicial complexes all whose vertex-links are in $\Sigma_k(d-1)$ (respectively in ${\mathcal S}_k(d-1)$). Thus,
members of ${\mathcal W}_k(d)$ (resp. ${\mathcal K}_k(d)$) are combinatorial manifolds (resp. homology manifolds)
without boundary. We have ${\mathcal W}_1(d) = {\mathcal K}_1(d)$ and ${\mathcal W}_k(d) \subseteq {\mathcal
K}_k(d)$ for $d\geq 2k$. The class ${\mathcal K}_k(d)$ is known as a generalized Walkup class (after D. W. Walkup
who considered the case $k=1$ in \cite{wa}). The class ${\mathcal W}_k(d)$ plays an important role in our recent
paper \cite{bd16} on tight triangulated manifolds. We also consider the class $\overline{{\mathcal W}}_k(d)$
(resp. $\overline{{\mathcal K}}_k(d)$) consisting of simplicial complexes all whose vertex-links are $k$-shelled
(resp. $k$-stacked homology) $(d - 1)$-balls. Thus members of $\overline{{\mathcal W}}_k(d)$ are combinatorial
manifolds with boundary while members of $\overline{{\mathcal K}}_k(d)$ are homology manifolds with boundary. We
have $\overline{{\mathcal W}}_1(d) = \overline{{\mathcal K}}_1(d)$ and $\overline{{\mathcal W}}_k(d) \subseteq
\overline{{\mathcal K}}_k(d)$ for $d \geq k+1$. Clearly the boundary of any member of $\overline{{\mathcal
W}}_k(d+1)$ (resp. $\overline{{\mathcal K}}_k(d+1)$) is in ${\mathcal W}_k(d)$ (resp. ${\mathcal K}_k(d)$) so
that we have the boundary map $\partial \colon \overline{{\mathcal W}}_k(d+1) \to {\mathcal W}_k(d)$ (resp.
$\partial \colon \overline{{\mathcal K}}_k(d+1) \to {\mathcal K}_k(d)$). Using the Murai-Nevo result quoted
above, we show that for $d\geq 2k+2$, the map $\partial \colon \overline{{\mathcal K}}_k(d+1) \to  {\mathcal
K}_k(d)$ is a bijection, and its inverse $M \mapsto \overline{M}$ has a simple combinatorial description.

In \cite{bd16}, we proved that, for a field $\mathbb{F}$, any $\mathbb{F}$-orientable $(k+1)$-neighborly member
of ${\mathcal W}_k(d)$ is $\mathbb{F}$-tight, provided $d\neq 2k+1$. Here we complement this result by proving
that an $\mathbb{F}$-orientable $(k+1)$-neighborly member of ${\mathcal W}_k(2k+1)$ is $\mathbb{F}$-tight if and
only if it is $k$-stacked.

In the final section of this paper, we present various examples, counterexamples and questions related to the
above results. For instance, we show that for each $k \geq 2$, there are $k$-stacked homology $d$-spheres which
are not even $(d+ 1)$-stellated (i.e., not combinatorial spheres) and $k$-stacked combinatorial $d$-spheres which
are not $d$-stellated. Recently, Klee and Novik \cite{kn} found an extremely beautiful construction of a
$(2d+4)$-vertex triangulation $M$ of $S^{\,k} \times S^{\,d-k}$ for all pairs $0\leq k\leq d$. We show that, for
$d\geq 2k$, these triangulations are in ${\mathcal W}_k(d)$. Klee and Novik obtained their triangulation $M$ as
the boundary complex of a triangulated $(d+ 1)$-manifold $\overline{M}$. For $d \geq 2k+2$, this is an instance
of our canonical construction $M \mapsto \overline{M}$. As an application, we show that, for $d\neq 2k$, the full
automorphism group of the Klee-Novik triangulation is a group of order $4d+8$, already found by these authors.
This makes it interesting to determine the full automorphism group of the Klee-Novik manifolds for $d=2k$.

\section{Bistellar moves and shelling moves}

A $d$-dimensional simplicial complex is called {\em pure} if all its maximal faces (called {\em facets}) are
$d$-dimensional. A $d$-dimensional pure simplicial complex is said to be a {\em weak pseudomanifold} if each of
its $(d - 1)$-faces is in at most two facets. For a $d$-dimensional weak pseudomanifold $X$, the {\em boundary}
$\partial X$ of $X$ is the pure subcomplex of $X$ whose facets are those $(d-1)$-dimensional faces of $X$ which
are contained in unique facets of $X$. The {\em dual graph} $\Lambda(X)$ of a weak pseudomanifold $X$ is the
graph whose vertices are the facets of $X$, where two facets are adjacent in $\Lambda(X)$ if they intersect in a
face of codimension one. A {\em pseudomanifold} is a weak pseudomanifold with a connected dual graph. A
$d$-dimensional weak pseudomanifold is called a {\em normal pseudomanifold} if each face of dimension $\leq d-2$
has a connected link. Since we include the empty set as a face, a normal pseudomanifold is necessarily connected.
All connected homology manifolds are automatically normal pseudomanifolds. We also know that every normal
pseudomanifold is a pseudomanifold (cf. \cite{bd9}).

For any two simplicial complexes $X$ and $Y$, their {\em join} $X \ast Y$ is the simplicial complex whose faces
are the disjoint unions of the faces of $X$ with the faces of $Y$. (Here we adopt the convention that the empty
set is a face of every simplicial complex.)

For a finite set $\alpha$, let $\overline{\alpha}$ (respectively $\partial \alpha$) denote the simplicial complex
whose faces are all the subsets (respectively, all proper subsets) of $\alpha$. Thus, if $\#(\alpha)= n \geq 2$,
$\overline{\alpha}$ is a copy of the standard triangulation $B^{\,n - 1}_n$ of the $(n- 1)$-dimensional ball, and
$\partial \alpha$ is a copy of the standard triangulation $S^{\,n-2}_n$ of the $(n-2)$-dimensional sphere. So,
for any two disjoint finite sets $\alpha$ and $\beta$, $\overline{\alpha}\ast\partial\beta$ and $\partial\alpha
\ast \overline{\beta}$ are two triangulations of a ball; they have identical boundaries, namely $(\partial
{\alpha}) \ast (\partial {\beta})$. We shall write $\dim(\alpha)$ for $\dim(\overline{\alpha}) = \#(\alpha) -1$.

A subcomplex $Y$ of a simplicial complex $X$ is said to be an {\em induced} (or {\em full}\,) subcomplex if every
face of $X$ contained in the vertex set of $Y$ is a face of $Y$. An induced subcomplex of $X$ with vertex set $U$
is denoted by $X[U]$. If $X$ is a $d$-dimensional simplicial complex with an induced subcomplex
$\overline{\alpha} \ast \partial \beta$ ($\alpha \neq \emptyset$, $\beta \neq \emptyset$) of dimension $d$ (thus,
$\dim(\alpha) + \dim(\beta) = d$), then $Y := (X\setminus (\overline{\alpha} \ast \partial \beta)) \cup
(\partial\alpha \ast \overline{\beta})$ is clearly another triangulation of the same topological space $|X|$. In
this case, $Y$ is said to be obtained from $X$ by the {\em bistellar move} $\alpha \mapsto \beta$. If
$\dim(\beta) = i$ ($0\leq i \leq d$), we say that $\alpha \mapsto \beta$ is a {\em bistellar move of index $i$}
(or an {\em $i$-move}, in short). Clearly, if $Y$ is obtained from $X$ by an $i$-move $\alpha \mapsto \beta$ then
$X$ is obtained from $Y$ by the (reverse) $(d- i)$-move $\beta \mapsto \alpha$. Notice that, in case $i=0$, i.e.,
when $\beta$ is a single vertex, we have $\partial \beta = \{\emptyset\}$ and hence $\overline{\alpha} \ast
\partial \beta = \overline{\alpha}$. Therefore, our requirement that $\overline{\alpha} \ast \partial \beta$ is
the induced subcomplex of $X$ on $\alpha \sqcup \beta$ means that $\beta$ is a new vertex, not in $X$. Thus, a
$0$-move creates a new vertex, and correspondingly a $d$-move deletes an old vertex. For $0 < i < d$, any
$i$-move preserves the vertex set; these are sometimes called the {\em proper bistellar moves}. For a thorough
treatment of bistellar moves, see \cite{bl}, for instance.

A triangulation $X$ of a manifold is called a {\em combinatorial manifold} if its geometric carrier $|X|$ is a
piecewise linear (pl) manifold with the pl structure induced from $X$. A combinatorial triangulation of a
sphere/ball  is called a {\em combinatorial sphere/ball} if it induces the standard pl structure (namely, that of
the standard sphere/ball) on its geometric carrier. Equivalently (cf. \cite{li, p}), a simplicial complex is a
combinatorial sphere (or ball) if it is obtained from a standard sphere (respectively, a standard ball) by a
finite sequence of bistellar moves. In general, a pure simplicial complex is a combinatorial manifold if and only
if the link of each of its vertices is a combinatorial sphere or combinatorial ball. (Recall that the {\em link}
of a vertex $x$ in a complex $X$, denoted by ${\rm lk}_X(x)$, is the subcomplex $\{\alpha \in X \, : \, x\not\in
\alpha, \alpha \sqcup\{x\} \in X\}$. Also, the {\em star} of $x$ in $X$, denoted by ${\rm st}_X(x)$, is the cone
$x \ast {\rm lk}_X(x)$. The {\em antistar} of $x$ in $X$, denoted by ${\rm ast}_X(x)$, is the subcomplex
$\{\alpha\in X \, : \, x\not\in\alpha\}$.) This leads us to introduce\,:

\begin{defn} \label{stellated-sphere}
{\rm For $0\leq k \leq d+1$, a $d$-dimensional simplicial complex $X$ is said to be {\em $k$-stellated} if $X$
may be obtained from $S^{\,d}_{d+2}$ by a finite sequence of bistellar moves, each of index $< k$. By convention,
$S^{\,d}_{d+2}$ is the only $0$-stellated simplicial complex of dimension $d$.

Clearly, for $0\leq k \leq l \leq d+1$, $k$-stellated implies $l$-stellated. All $k$-stellated simplicial
complexes are combinatorial spheres. We let $\Sigma_k(d)$ denote the class of all $k$-stellated $d$-spheres. By
Pachner's theorem (\cite{p}), $\Sigma_{d+1}(d)$ consists of all combinatorial $d$-spheres.
 }
\end{defn}

By definition, $X\in \Sigma_k(d)$ if and only if there is a sequence $X_0, X_1, \dots, X_n$ of $d$-dimensional
simplicial complexes such that $X_0 = S^{\,d}_{d + 2}$, $X_n = X$ and, for $0 \leq j < n$, $X_{j + 1}$ is
obtained from $X_j$ by a single bistellar move of index $\leq k-1$. The smallest such integer $n$ is said to be
the {\em length} of $X\in \Sigma_k(d)$ and is denoted by $l(X)$. For $X, Y \in \Sigma_k(d)$, we say that $Y$ is
shorter than $X$ if $l(Y) < l(X)$. Thus, $S^{\,d}_{d + 2}$ is the unique shortest member of $\Sigma_k(d)$ (of
length 0), and every other member of $\Sigma_k(d)$ can be obtained from a shorter member by a single bistellar
move of index $< k$. Thus, induction on the length is a natural method for proving results about the class
$\Sigma_k(d)$.

Let $X$, $Y$ be two pure simplicial complexes of dimension $d$. We say that $X$ is obtained from $Y$ by the {\em
shelling move} $\alpha \leadsto \beta$ if $\alpha$ and $\beta \neq \emptyset$ are disjoint faces of $X$ such that
(i) $Y \subseteq X$, and $\alpha \sqcup \beta$ is the only facet of $X$ which is not a facet of $Y$, and (ii) the
induced subcomplex of $Y$ on the vertex set of $\alpha \sqcup \beta$ is $\overline{\alpha}\ast \partial \beta$.
If $\dim(\beta) = i$, we say that the shelling move $\alpha \leadsto \beta$ is of index $i$. (Clearly,
$\dim(\alpha) + \dim(\beta) = d - 1$, so that $0 \leq i \leq d$).

We say that a $d$-dimensional simplicial complex $X$ {\em shellable} if $X$ is obtained from the standard
$d$-ball $B^{\,d}_{d + 1}$ by a finite sequence of shelling moves. Clearly, each shelling move increases the
number of facets by one, so that - when $X$ is shellable, the number of shelling moves needed to obtain $X$ from
$B^{\,d}_{d + 1}$ is one less than the number of facets of $X$.

Let $X$ and $Y$ be $d$-dimensional pseudomanifolds. If $X$ is obtained from $Y$ by the shelling move $\alpha
\leadsto \beta$ then $X = Y \cup \overline{\alpha \sqcup \beta}$, $Y \cap \overline{\alpha \sqcup \beta} =
\overline{\alpha}\ast \partial \beta$. (Since $X$ is a pseudomanifold, it follows that $\overline{\alpha}\ast
\partial \beta \subseteq \partial Y$.) If the move is of index $<d$, then $\overline{\alpha} \ast \partial \beta$
is a combinatorial $(d-1)$-ball; if it is of index $d$ (so that $\alpha = \emptyset$), $\overline{\alpha} \ast
\partial \beta$ ($=\partial \beta$) is a combinatorial $(d-1)$-sphere. Therefore, if $Y$ is a combinatorial
$d$-ball, then $X$ is also a combinatorial $d$-ball in case the shelling move is of index $<d$, and $X$ is a
combinatorial $d$-sphere if the shelling move is of index $d$. (Also note that $Y$ can't be a combinatorial
sphere since a $d$-dimensional pseudomanifold without boundary can't be properly contained in a
$d$-pseudomanifold with or without boundary.) From these observations, it is immediate by an induction on the
number of facets that a shellable pseudomanifold is either a combinatorial ball or a combinatorial sphere. (This
result appears to be due to Danaraj and Klee \cite{dk}.) Also if $X$ is a shellable $d$-pseudomanifold, then
among the shelling moves used to obtain $X$ from $B^{\,d}_{d +1}$, only the last move can be of index $d$; this
happens if and only if $X$ is a $d$-sphere. These considerations lead us to introduce\,:

\begin{defn} \label{k-shelled-ball}
{\rm For $0\leq k \leq d$, a $d$-dimensional pseudomanifold is said to be {\em $k$-shelled} if it may be obtained
from the standard $d$-ball $B^{\,d}_{d+1}$ by a finite sequence of shelling moves, each of index $< k$. By
convention, $B^{\,d}_{d+1}$ is the only $0$-shelled pseudomanifold of dimension $d$. }
\end{defn}

Clearly, all $k$-shelled pseudomanifolds are combinatorial balls. Also, for $0\leq k \leq l\leq d$, $k$-shelled
implies $l$-shelled. By $\widehat{\Sigma}_k(d)$, $0\leq k\leq d$, we denote the class of all $k$-shelled
$d$-balls. Thus $\widehat{\Sigma}_d(d)$ consists of all the shellable $d$-balls. Note that, while all shellable
balls are combinatorial balls, the converse is false.

Unlike the case of bistellar moves, the reverse of a shelling move is not a shelling move. Nonetheless, the two
notions are closely related, as the following lemma shows. (This lemma seems to be well known to experts; but we
could not find a reference.)

\begin{lemma}\label{lemma-1}
If a homology $(d+1)$-ball $X$ is obtained from a homology $(d+1)$-ball $Y$ by a shelling move $\alpha \leadsto
\beta$ of index $i \leq d$ then the homology $d$-sphere $\partial X$ is obtained from the homology $d$-sphere
$\partial Y$ by the bistellar move $\alpha \mapsto \beta$ of index $i$.
\end{lemma}

\noindent {\bf Proof.} Let $\sigma = \alpha \sqcup \beta$. Thus, $\sigma$ is the only facet of $X$ which is not
in $Y$. Since $Y \subseteq X$ are $(d+1)$-dimensional pseudomanifolds, it follows that (i) a boundary $d$-face of
$Y$ is not a boundary $d$-face of $X$ if and only if (it is a face of $Y$ and) it is contained in $\sigma$, i.e.,
if and only if it is a facet of $\overline{\alpha} \ast \partial \beta$, and (ii) a boundary $d$-face of $X$ is
not a face of $Y$ if and only if it is a facet of $\overline{\beta}\ast
\partial \alpha$. Since $\partial X$ and $\partial Y$ are pure
simplicial complexes of dimension $d$, the result follows. \hfill $\Box$

\medskip

As an immediate consequence of this lemma, we have\,:

\begin{cor} \label{coro-1}
If $B$ is a $k$-shelled $(d+1)$-ball then $\partial B$ is a $k$-stellated $d$-sphere.
\end{cor}

For a simplicial complex $X$, say of dimension $d$, and a non-negative integer $m \leq d$, the {\em $m$-skeleton}
of $X$, denoted by ${\rm skel}_m(X)$, is the subcomplex of $X$ consisting of all its faces of dimension $\leq m$.
We recall\,:

\begin{defn} \label{k-stacked-ball}
{\rm For $0\leq k \leq d+1$, a homology $(d+1)$-dimensional ball $B$ is said to be {\em $k$-stacked} if all the
faces of $B$ of codimension (at least) $k+1$ lie in its boundary; i.e., if ${\rm skel}_{d - k}(B) = {\rm skel}_{d
- k}(\partial B)$. A homology $d$-sphere $S$ is said to be {\em $k$-stacked} if there is a $k$-stacked homology
$(d+ 1)$-ball $B$ such that $\partial B = S$. We let ${\mathcal S}_k(d)$ and  $\widehat{\mathcal S}_k(d)$ denote
the class of all $k$-stacked homology $d$-spheres and of all $k$-stacked homology $d$-balls respectively. }
\end{defn}

Clearly, we have ${\mathcal S}_0(d) \subseteq {\mathcal S}_1(d) \subseteq \cdots \subseteq {\mathcal S}_d(d)$ and
$\widehat{\mathcal S}_0(d) \subseteq \widehat{\mathcal S}_1(d) \subseteq \cdots \subseteq \widehat{\mathcal
S}_d(d)$. Trivially, the standard $d$-ball is the only member of $\widehat{\mathcal S}_0(d)$, and hence the
standard $d$-sphere is the only member of ${\mathcal S}_0(d)$. Our first Theorem shows that ${\mathcal S}_d(d)$
consists of all the homology $d$-spheres. Notice that, trivially, $\widehat{\mathcal S}_d(d)$ consists of all
homology $d$-balls.

\begin{theo} \label{theo-1}
Every homology $d$-sphere is $d$-stacked.
\end{theo}

\noindent {\bf Proof.} Let $S$ be a homology $d$-sphere. Fix a vertex $x$ of $S$. Let $A_x$ be the antistar of
$x$ in $S$. Set $B_x = \overline{\{x\}} \ast A_x$. It is shown in Lemma 9.1 of \cite{bd9} that $B_x$ is a
homology $(d+1)$-ball. Clearly, $B_x$ has the same vertex set as $S = \partial B_x$. Therefore, $B_x$ is a
$d$-stacked homology $(d+1)$-ball and (hence) $S$ is a $d$-stacked homology $d$-sphere. \hfill $\Box$

\begin{theo} \label{theo-2}
Let $B$ be a homology $(d+1)$-ball. Then $B$ is $k$-shelled if and only if $B$ is shellable and $k$-stacked.
\end{theo}

\noindent {\bf Proof.} Suppose $B$ is $k$-shelled. Then, of course, $B$ is shellable. We prove that $B$ is
$k$-stacked by induction on the number of facets of $B$. If $B$ has only one facet then $B = B^{\,d + 1}_{d +
2}$, the standard ball, and the result is trivial. Otherwise, $B$ is obtained from a $k$-shelled ball
$B^{\hspace{.2mm} \prime}$ (with one less facet) by a single shelling move $\alpha \leadsto \beta$ of index $\leq
k - 1$. By induction hypothesis, ${\rm skel}_{d - k}(B^{\hspace{.2mm} \prime}) = {\rm skel}_{d - k}(\partial
B^{\hspace{.2mm} \prime})$, and by Lemma \ref{lemma-1}, $\partial B$ is obtained from $\partial B^{\hspace{.2mm}
\prime}$ by the bistellar move $\alpha \mapsto \beta$ of index $\leq k-1$.

Let $\gamma$ be a face of $B$ of dimension $\leq d-k$. Since $\dim(\alpha) \geq d-k+1$, $\gamma \not \supseteq
\alpha$. If $\gamma$ is a face of $B^{\hspace{.2mm} \prime}$ then (as $B^{\hspace{.2mm} \prime}$ is $k$-stacked),
$\gamma \in \partial B^{\hspace{.2mm} \prime}$. Since $\gamma \not\supseteq \alpha$, and $\partial B$ is obtained
from $\partial B^{\hspace{.2mm} \prime}$ by the bistellar move $\alpha \mapsto \beta$, it follows that $\gamma
\in \partial B$. If, on the other hand, $\gamma$ is not a face of $B^{\hspace{.2mm} \prime}$ then $\beta
\subseteq \gamma \subseteq \alpha\sqcup\beta$ and hence we have $\gamma \in \overline{\beta}\ast \partial \alpha
\subseteq \partial B$. Thus $\gamma \in \partial B$ in either case. So, $B$ is $k$-stacked. This proves the
``only if\," part.

The ``if part\," is also proved by induction on the number of facets of $B$. Suppose $B$ is a $k$-stacked
shellable $(d+1)$-ball. If $B = B^{\,d+1}_{d+2}$, then $B$ is vacuously $k$-shelled. Else, $B$ is obtained from a
shellable $(d+1)$-ball $B^{\hspace{.2mm} \prime}$ (with one less facet) by a single shelling move $\alpha
\leadsto \beta$. By Lemma \ref{lemma-1}, $\partial B$ is obtained from $\partial B^{\hspace{.2mm} \prime}$ by the
bistellar move $\alpha \mapsto \beta$. Hence $\alpha \not\in \partial B$ but $\alpha \in B$. Since $B$ is
$k$-stacked, it follows that $\dim(\alpha) \geq d-k+1$, and hence $\dim(\beta) \leq k-1$. Thus, the shelling move
$\alpha \leadsto \beta$ is of index $\leq k-1$. Let $\gamma \in B^{\hspace{.2mm} \prime}$, $\dim(\gamma) \leq
d-k$. Since $B^{\hspace{.2mm} \prime}\subseteq B$, it follows that $\gamma \in B$. Since $\dim(\gamma) \leq d-k$
and $B$ is $k$-stacked, it follows that $\gamma \in \partial B$. As $\beta \not\in B^{\hspace{.2mm} \prime}$ and
$\gamma \in B^{\hspace{.2mm} \prime}$, we also have $\gamma \not\supseteq \beta$. Thus $\gamma \not\supseteq
\beta$, $\gamma \in \partial B$ and $\partial B$ is obtained from $\partial B^{\hspace{.2mm} \prime}$ by the
bistellar move $\alpha \mapsto \beta$. Hence $\gamma \in \partial B^{\hspace{.2mm} \prime}$. This shows that
$B^{\hspace{.2mm} \prime}$ is $k$-stacked. As $B^{\hspace{.2mm} \prime}$ is $k$-stacked and shellable, the
induction hypothesis implies that $B^{\hspace{.2mm} \prime}$ is $k$-shelled. Since $B$ is obtained from
$B^{\hspace{.2mm} \prime}$ by a shelling move of index $\leq k-1$, it follows that $B$ is also $k$-shelled. This
completes the induction. \hfill $\Box$

\bigskip

Thus we have $\widehat{\Sigma}_k(d) \subseteq \widehat{\mathcal S}_k(d)$. Our next result gives a one-sided
relationship between $k$-stacked spheres and $k$-stacked balls on one hand, and between $k$-stellated spheres and
$k$-shelled balls on the other hand.

\begin{theo} \label{theo-3}
Let $B$ be a homology ball. \vspace{-2mm}
\begin{enumerate}
\item[$(a)$] If $B$ is $k$-stacked then there is a $k$-stacked homology sphere $S$ such that $B$ is the antistar
of a vertex in $S$. \vspace{-2mm}

\item[$(b)$] If $B$ is $k$-shelled then there is a $k$-stellated sphere $S$ such that $B$ is the antistar of a
vertex in $S$.
\end{enumerate}
\end{theo}

\noindent {\bf Proof.} Let $x$ be a new vertex (not in $B$), and set $S := B \cup (x\ast \partial B)$. (Notice
that, since $S$ is to be a $d$-pseudomanifold without boundary and $B$ is a $d$-pseudomanifold with boundary,
this is the only choice of $S$ so that $B$ is the antistar of a vertex $x$ in $S$.) Clearly, $S = \partial B_0$,
where $B_0 = x \ast B$. Therefore, to prove the result, it is enough to show that if $B$ is $k$-stacked
(respectively $k$-shelled) then so is $B_0$. But, this is trivial. \hfill $\Box$

\bigskip

Next we present a characterization of $k$-stellated spheres of dimension $\geq 2k-1$.

\begin{theo} \label{theo-4}
A homology sphere of dimension $\geq 2k -1$ is $k$-stellated if and only if it is the boundary of a $k$-shelled
ball. In consequence, all $k$-stellated spheres of dimension $\geq 2k-1$ are $k$-stacked.
\end{theo}

\noindent {\bf Proof.} The ``if\," part is Corollary \ref{coro-1} (which holds in all dimensions). We prove the
``only if\," part by induction on the length $l(S)$ of a $k$-stellated sphere $S$ of dimension $d\geq 2k-1$. If
$l(S) = 0$ then $S = S^{\,d}_{d+2}$ is the boundary of $B^{\,d+1}_{d+2}$. So, let $l(S) > 0$. Then $S$ is
obtained from a shorter member $S^{\hspace{.2mm} \prime}$ of $\Sigma_k(d)$ by a single bistellar move $\alpha
\mapsto \beta$ of index $\leq k-1$. By induction hypothesis, there is a $k$-shelled $(d+1)$-ball
$B^{\hspace{.2mm}\prime}$ such that $\partial B^{\hspace{.2mm}\prime} = S^{\hspace{.2mm}\prime}$. The induced
subcomplex of $S^{\hspace{.2mm} \prime}$ on the vertex set $\alpha\sqcup \beta$ is $\overline{\alpha} \ast
\partial \beta \subseteq S^{\hspace{.2mm} \prime} \subseteq B^{\hspace{.2mm} \prime}$. Since $\dim(\beta) \leq
k-1 \leq d-k$, $\beta \not\in S^{\hspace{.2mm} \prime} = \partial B^{\hspace{.2mm} \prime}$ and (by Theorem
\ref{theo-2}) $B^{\hspace{.2mm} \prime}$ is $k$-stacked, it follows that $\beta \not\in B^{\hspace{.2mm}
\prime}$. Thus, the induced subcomplex of $B^{\hspace{.2mm} \prime}$ on $\alpha \sqcup \beta$ is also
$\overline{\alpha} \ast \partial \beta$. So, $B^{\hspace{.2mm} \prime}$ admits the shelling move $\alpha \leadsto
\beta$ of index $\leq k-1$. Let $B$ be the $(d+1)$-ball obtained from $B^{\hspace{.2mm} \prime}$ by this move.
Since $B^{\hspace{.2mm} \prime}$ is $k$-shelled, so is $B$. By Lemma \ref{lemma-1}, $\partial B$ is obtained from
$S^{\hspace{.2mm} \prime} = \partial B^{\hspace{.2mm} \prime}$ by the bistellar move $\alpha \mapsto \beta$. That
is, $\partial B = S$. This completes the induction. The second statement is now immediate from the first
statement and Theorem \ref{theo-2}. \hfill $\Box$

\bigskip

Recall that a triangulated $d$-sphere is said to be {\em polytopal} if it is isomorphic to the boundary complex
of a $(d+1)$-dimensional simplicial polytope. A simplicial complex $X$ is said to be {\em $l$-neighbourly} if any
$l$ vertices of $X$ form a face of $X$.

\begin{theo} \label{theo-5}
Let $S$ be a $(k+1)$-neighbourly polytopal sphere of dimension $d$. Then $S$ is $(d-k)$-stellated.
\end{theo}

\noindent {\bf Proof.} Fix a vertex $x$ of $S$. Let $A$ be the antistar of $x$ in $S$. By Bruggesser-Mani (cf.
\cite[Theorem 8.12]{z95}) $A$ is a shellable $d$-ball. Hence, $x \ast A$ is a shellable $(d+1)$-ball. Clearly,
$\partial(x\ast A) = S$. Since $S$ is $(k+1)$-neighbourly, $x \ast A$ is $(d-k)$-stacked. Hence, by Theorem
\ref{theo-2}, $x \ast A$ is $(d-k)$-shelled. Therefore, by Corollary \ref{coro-1}, $S$ is $(d-k)$-stellated.
\hfill $\Box$

\bigskip

Recall that a {\em missing face} of dimension $l$ in a simplicial complex $X$ is a set $\alpha$ consisting of
$l+1$ vertices of $X$ such that $\alpha$ is not a face of $X$, but all proper subsets of $\alpha$ are faces of
$X$. In other words, $\alpha$ is a missing $l$-face of $X$ if and only if the induced subcomplex $X[\alpha]$
(with vertex set $\alpha$) of $X$ is a standard sphere $S^{\,l - 1}_{l+1}$. In \cite{bd13v2}, we had proved the
special case of the following result for polytopal balls. Also see Corollary 3.2 in \cite{mn2}.

The proof of the following result closely follows that of Theorem 2.3 (ii) of \cite{mn}. So, we include a brief
sketch of the proof.

\begin{lemma} \label{lemma-2}
Let $B$ be a $k$-stacked homology ball. Then all the missing faces of $B$ have dimension $\leq k$.
\end{lemma}

\noindent {\bf Proof.} Let $\dim(B) = d+1$, and $S = \partial B$. Thus, $S$ is a homology $d$-sphere with ${\rm
skel}_{d-k}(S) = {\rm skel}_{d-k}(B)$. Take a new vertex $x$ and form the cone $\widehat{B} = x\ast B$. Let's put
$\widehat{S} = \partial \widehat{B} = B \cup (x \ast S)$. Let $V$ be the vertex set of $S$ and $\widehat{V} = V
\sqcup \{x\}$ be the vertex set of $\widehat{S}$.

In the following, we fix a field $\mathbb{F}$ such that $S$ (and hence also $\widehat{S}$) is an
$\mathbb{F}$-homology sphere. All homologies used below are simplicial homologies with coefficients in
$\mathbb{F}$.

Let $\alpha \subseteq V$ with $\#(\alpha) = l+1$ (say), where $l \geq k+1$.  Let $\beta = \widehat{V}
\setminus\alpha$ and $\gamma = V\setminus \alpha$. Thus $\beta = \gamma \sqcup \{x\}$. Since $d + 1 - l \leq
d-k$, we have ${\rm skel}_{d-l+1}(B) = {\rm skel}_{d - l + 1}(S)$. Hence ${\rm skel}_{d-l+1}(\widehat{S}) = {\rm
skel}_{d -l+1}(x\ast S)$. It follows that ${\rm skel}_{d-l+ 1}(\widehat{S}[\beta]) = {\rm skel}_{d-l+1}(x\ast
S[\gamma])$. Since $x \ast S[\gamma]$ is a cone and $x\ast S \subseteq \widehat{S}$, this implies
$Z_{d-l+1}(x\ast S[\gamma]) = B_{d-l+1}(x\ast S[\gamma]) \subseteq B_{d-l + 1}(\widehat{S}[\beta]) \subseteq
Z_{d-l+ 1}(\widehat{S}[\beta]) = Z_{d-l+1}(x\ast S[\gamma])$. This implies $H_{d-l+ 1}(\widehat{S}[\beta])
=\{0\}$.

Since $\widehat{S}$ is an $\mathbb{F}$-homology $(d+1)$-sphere and $\beta$ is the complement of $\alpha$ in the
vertex set of $\widehat{S}$, simplicial Alexander duality (see, for example, Lemma 4.1 in \cite{bd8}) and the
exact sequence for pairs imply that $H_{l-1}(\widehat{S}[\alpha]) = H_{d-l+1}(\widehat{S}[\beta]) = \{0\}$. But,
$\alpha \subseteq V$ and $B = \widehat{S}[V]$. Therefore, we have $B[\alpha] = \widehat{S}[\alpha]$. Thus, we get
$H_{l-1}(B[\alpha]) = \{0\}\neq H_{l-1}(S^{l-1}_{l+1})$. Hence $B[\alpha] \neq S^{l-1}_{l+1}$. Thus, $\alpha$ is
not a missing face of $B$. \hfill $\Box$

\bigskip

The following result is essentially Theorem 2.3 (ii) of Murai and Nevo \cite{mn}.

\bigskip

\noindent {\bf Notation\,:} For a set $\alpha$ and a non-negative integer $m$, $\binom{\alpha}{\leq \, m\,}$ will
denote the collection of all subsets of $\alpha$ of size $\leq m$.

\begin{theo} \label{theo-6}
Let $S$ be a $k$-stacked homology sphere of dimension $d \geq 2k$, say with vertex set $V$. Then there is a
unique $k$-stacked homology $(d+ 1)$-ball $\overline{S}$ whose boundary is $S$. $($If further, $S$ is a
$k$-stellated sphere then, by Theorems $\ref{theo-2}$ and $\ref{theo-4}$, $\overline{S}$ is actually
$k$-shelled.$)$ It is given by the formula \vspace{-2mm}
\begin{eqnarray} \label{eq1}
\overline{S} = \left\{\alpha \subseteq V \, : \, \binom{\alpha}{\leq \, k+1\,} \subseteq S\right\}.
\end{eqnarray}
\end{theo}

\noindent {\bf Proof.} Let $B$ be a homology $(d+1)$-ball such that $\partial B=S$ and ${\rm skel}_{d-k}(B) =
{\rm skel}_{d- k}(S)$.  We must show that $B = \overline{S}$. Since $d \geq 2k$, we have ${\rm skel}_k(B) = {\rm
skel}_{k}(S) \subseteq S$, and therefore, by the definition of $\overline{S}$, we have $B \subseteq
\overline{S}$. If $B \neq \overline{S}$, then choose an inclusion minimal member $\alpha$ of
$\overline{S}\setminus B$. Then $\alpha$ is a missing face of $B$. Therefore, by Lemma \ref{lemma-2},
$\dim(\alpha) \leq k$. Then $\alpha \in {\rm skel}_{k}(\overline{S}) \subseteq S \subseteq B$. Thus $\alpha\in
B$; a contradiction. \hfill $\Box$

\bigskip

In \cite{bd13v2}, we had proved two special cases of Theorem \ref{theo-6}\,: for $k$-stellated spheres and for
$k$-stacked polytopal spheres. In Proposition 3.6 of \cite{ka}, Kalai proved the special case of the following
corollary for polytopal spheres. Also, in Corollary 4.8 of \cite{na}, Nagel proved the special case of this
corollary for homology spheres with the Weak Lefschetz Property (WLP). Conjecturally, all homology spheres have
WLP. However, our proof is unconditional. (This is also proved in Remark 4.5 of \cite{mn2}.)

\begin{cor} \label{coro-2}
For $k \leq e \leq d-k-1$, a $k$-stacked homology $d$-sphere does not have any standard $e$-sphere as an induced
subcomplex. In consequence, such a $d$-sphere does not admit any bistellar move of index $i$ for $k+1 \leq i \leq
d-k$.
\end{cor}

\noindent {\bf Proof.} Notice that a homology sphere $S$ admits a bistellar move $\alpha \mapsto \beta$ of index
$i$ if and only if it has $\overline{\alpha}\ast \partial\beta$ as an induced subcomplex. In this case, it has
the standard $(i-1)$-sphere $\partial \beta$ as the induced subcomplex on $\beta$. So, the second statement is
immediate from the first. The first statement is vacuously true unless $d \geq 2k+1$. So, to prove it, we may
assume $d \geq 2k+1$. By Theorem \ref{theo-6}, we have ${\rm skel}_{d- k}({S}) = {\rm skel}_{d-k}(\overline{S})$.
Hence any induced standard sphere of dimension $e \leq d-k-1$ in $S$ is also an induced standard sphere of
$\overline{S}$, so that $e \leq k-1$ by Lemma \ref{lemma-2}. This proves the first statement. \hfill $\Box$

\bigskip

If $S$ is a $k$-stellated $d$-sphere, other than the standard sphere, then $S$ is obtained from a shorter
$k$-stellated $d$-sphere by a bistellar move of index $\leq k-1$. Hence such a sphere admits the reverse move,
which is a bistellar move of index $\geq d-k+1$. In consequence, such a sphere always has an induced subcomplex
isomorphic to a standard sphere of some dimension $\geq d-k$. In this sense, Corollary \ref{coro-2} is best
possible. Indeed, it is easy to prove by induction on the length that if $d \geq 2k-2$ and $S$ is a $k$-stellated
$d$-sphere which is not $(k-1)$-stellated, then $S$ has an $S^{\,d-k}_{d-k+2}$ as an induced subcomplex.

In the following proof (and also later) we use the notation $V(X)$ for the vertex set of a simplicial complex
$X$.

\begin{theo} \label{theo-7}
For a normal pseudomanifold $X$, the following are equivalent\,: \vspace{-1mm}
\begin{enumerate}
\item[$(i)$] $X$ is a $1$-shelled ball, \vspace{-1mm} \item[$(ii)$] $X$ is a $1$-stacked ball, \vspace{-1mm}
\item[$(iii)$] $X$ is a $1$-stacked $R$-homology ball for some commutative ring $R$, \vspace{-1mm} \item[$(iv)$]
$\Lambda(X)$ is a tree.
\end{enumerate}
\end{theo}

\noindent {\bf Proof.} Let $X$ be of dimension $d+1\geq 1$.

\smallskip

\noindent $(i) \Rightarrow (ii)$\,: Follows from Theorem \ref{theo-2}.

\smallskip

\noindent $(ii) \Rightarrow (iii)$\,: Follows from the fact that triangulated balls are homology balls.

\smallskip

\noindent $(iii) \Rightarrow (iv)$\,: The result is trivial for dimension 1. So, assume that $d+ 1 \geq 2$. If
$X$ has only one facet then the result is trivial. So, assume that $X$ is a 1-stacked $R$-homology ball with at
least two facets. Since $X$ is a homology ball, $\Lambda(X)$ is connected. To prove that $\Lambda(X)$ is a tree,
it suffices to show that each edge of $\Lambda(X)$ is a cut edge (i.e., deletion of any edge from $\Lambda(X)$
disconnects the graph). Let $e_0 = \sigma_1\sigma_2$ be an edge of $\Lambda(X)$. Then $\gamma := \sigma_1\cap
\sigma_2$ is an interior $d$-face of $X$; i.e., $\gamma \not\in S := \partial X$. Since ${\rm skel}_{d-1}(X) =
{\rm skel}_{d-1}(S)$, $\partial\gamma\subseteq S$. Thus, $\partial \gamma$ is an induced $S^{\,d-1}_{d+1}$ in the
$d$-sphere $S$. By Lemma 3.3 of \cite{bd9}, $S$ is obtained from a $d$-dimensional weak pseudomanifold
$\widetilde{S}$ (without boundary) by an elementary handle addition.

\smallskip

\noindent {\em Claim 1.} $\widetilde{S}$ is disconnected.

\smallskip

Let $S$ be obtained from $\widetilde{S}\setminus\{\alpha_1, \alpha_2\}$ by identifying vertices of simplices
$\alpha_1, \alpha_2$, where $\alpha_1$, $\alpha_2$ are disjoint facets in $\widetilde{S}$ (see \cite{bd9}). If
$\widetilde{S}$ is connected then, by using the exact sequence of pairs, we get $R \cong H_1(\widetilde{S},
\alpha_1\cup\alpha_2; R) = H_1(\widetilde{S}/\alpha_1 \cup\alpha_2; R) = H_1(S/\partial\gamma; R) = H_1(S,
\partial\gamma; R)$. This implies that $H_1(S; R) \cong R$, a contradiction. This proves the claim.

Since $\widetilde{S}$ is disconnected, by Lemma 3.3 of \cite{bd9}, $\widetilde{S}$ has exactly two components,
say $S_1$ and $S_2$. Then $S = S_1 \# S_2$ (connected sum) and $V(S_1) \cap V(S_2) = \gamma$. For $1\leq i \leq
2$, let $U_i$ be the set of facets of $X$ contained in $V(S_i)$. Since $V(S_1) \cap V(S_2) = \gamma$, it follows
that $U_1 \cap U_2 = \emptyset$.

If the dimension $d+ 1 = 2$ then $\gamma$ is an edge and it clearly divides the 2-disc $X$ into two parts and the
triangles (facets) in one part are in $U_1$ and the triangles in the other part are in $U_2$. Now, assume that
$d+1 \geq 3$. Let $uv$ be an edge of $X$. Since $d+1 \geq 3$, $uv \in S$ and hence (since $S = S_1\# S_2$) $uv
\in S_1$ or $uv \in S_2$. Therefore, $u, v\in V(S_1)$ or $u, v\in V(S_2)$. This implies that for any facet
$\sigma$ in $X$, either all the vertices of $\sigma$ are in $V(S_1)$ or all the vertices of $\sigma$ are in
$V(S_2)$. Thus, any facet in $X$ is in $U_1$ or in $U_2$. Thus (for any dimension $d+1 \geq 2$), $U_1 \sqcup U_2$
is a partition of the vertex set of the dual graph $\Lambda(X)$. Any facet $\sigma$ of $X$ containing a $d$-face
$\alpha \neq \gamma$ of $S_1$ is in $U_1$. So, $U_1 \neq \emptyset$. Similarly, $U_2 \neq \emptyset$.

Now, let $e = \alpha_1 \alpha_2$ be an edge of $\Lambda(X)$ with $\alpha_i\in U_i$, $i = 1, 2$. Then $\alpha :=
\alpha_1\cap \alpha_2 \subseteq V(S_i)$ for $i=1, 2$. Hence $\alpha \subseteq V(S_1) \cap V(S_2) = \gamma$ and
therefore $\alpha = \gamma$. So, $e = e_0$. Thus, $e_0$ is the unique edge of $\Lambda(X)$ with one end in $U_1$
and other end in $U_2$. So, $e_0$ is a cut edge of $\Lambda(X)$. Since $e_0$ was an arbitrary edge of
$\Lambda(X)$, this proves that $\Lambda(X)$ is a tree.

\smallskip

\noindent $(iv) \Rightarrow (i)$\,: Suppose $\Lambda(X)$ is a tree. We prove that $X$ is 1-shelled by induction
on the number of facets of $X$ (i.e., the number of vertices of $\Lambda(X)$). This is trivial if $X$ has only
one facet, i.e., $X = B^{\,d+1}_{d+2}$. So, assume $\Lambda(X)$ is a tree with at least two vertices. Then
$\Lambda(X)$ has a vertex $\sigma$ of degree 1 (leaf). Let $\sigma^{\hspace{.2mm}\prime}$ be the unique neighbour
of $\sigma$ in $\Lambda(X)$. Let $X^{\hspace{.2mm}\prime}$ be the pure simplicial complex whose facets are those
of $X$ other than $\sigma$.

\smallskip

\noindent {\em Claim 2.} If $\gamma = \sigma \cap \sigma^{\hspace{.2mm}\prime}$ then $X^{\hspace{.2mm}\prime}
\cap \overline{\sigma} = \overline{\gamma}$.

\smallskip

Let $\sigma = \gamma \sqcup \{u\}$. To prove Claim 2, it is sufficient to show that $u \not\in X^{\hspace{.2mm}
\prime}$. If possible let $u\in X^{\hspace{.2mm}\prime}$. Let $\alpha\subseteq X^{\hspace{.2mm}\prime} \cap
\overline{\sigma}$ be a maximal simplex containing $u$. Since $\sigma$ is a leaf in $\Lambda(X)$, $\dim(\alpha)
\leq d-2$. Clearly, ${\rm lk}_X(\alpha) = {\rm lk}_{X^{\hspace{.2mm}\prime}}(\alpha) \sqcup \overline{\sigma
\setminus \alpha}$ and ${\rm lk}_{X^{\hspace{.2mm}\prime}}(\alpha) \cap \overline{\sigma\setminus \alpha} =
\emptyset$. This is a contradiction since $X$ is a normal pseudomanifold. This proves the claim.

Clearly, $X^{\hspace{.2mm} \prime}$ is a normal pseudomanifold and $\Lambda(X^{\hspace{.2mm}\prime})$ is the tree
obtained from the tree $\Lambda(X)$ by deleting the end vertex $\sigma$ and the edge
$\sigma\sigma^{\hspace{.2mm}\prime}$. Therefore, by induction hypothesis, $X^{\hspace{.2mm}\prime}$ is a
1-shelled ball. By Claim 2, $X$ is obtained from $X^{\hspace{.2mm}\prime}$ by the shelling move $\gamma \leadsto
\{u\}$ of index 0. Therefore, $X$ is also a 1-shelled ball. \hfill $\Box$

\bigskip

Thus a homology ball is 1-stacked if and only if it is 1-shelled. So, $\widehat{\Sigma}_1(d) =\widehat{\mathcal
S}_1(d)$.  Now, Theorems \ref{theo-4} and \ref{theo-7} imply\,:

\begin{cor}  \label{coro-3}
A homology sphere is $1$-stellated if and only if it is $1$-stacked.
\end{cor}

The `only if' part of this corollary follows from the definitions and the `if' is also a consequence of Corollary
8.4 of \cite{ka87}.

Next we introduce\,:

\begin{defn} \label{defn-Wkd}
{\rm For $0 \leq k \leq d$, ${\mathcal W}_k(d)$ consists of the connected simplicial complexes of dimension $d$
all whose vertex-links are $k$-stellated $(d-1)$-spheres, and ${\mathcal K}_k(d)$ consists of the connected
simplicial complexes of dimension $d$ all whose vertex-links are $k$-stacked $(d-1)$-spheres. }
\end{defn}

Thus, members of ${\mathcal W}_k(d)$ are combinatorial manifolds; the members of ${\mathcal K}_k(d)$ are homology
manifolds. In consequence of Corollary \ref{coro-3}, we have\,:

\begin{cor}  \label{coro-4}
${\mathcal W}_1(d) ={\mathcal K}_1(d)$.
\end{cor}

In consequence of Theorem \ref{theo-4}, we have\,:

\begin{cor} \label{coro-5}
${\mathcal W}_k(d) \subseteq {\mathcal K}_k(d)$ for $d \geq 2k$.
\end{cor}

\begin{theo} \label{theo-8}
$(a)$ All $k$-stellated $d$-spheres belong to the class ${\mathcal W}_k(d)$. $(b)$ All $k$-stacked homology
$d$-spheres belong to the class ${\mathcal K}_k(d)$.
\end{theo}

\noindent {\bf Proof.} Let $S$ be a $k$-stellated $d$-sphere. We need to show that all the vertex-links of $S$
are $k$-stellated. Again, the proof is by induction on the length $l(S)$ of $S$. If $l(S) = 0$ then $S =
S^{\,d}_{d +2}$, and all its vertex links are $S^{\,d-1}_{d +1}$, so we are done. Therefore, let $l(S) >0$. Then
$S$ is obtained from a shorter $k$-stellated $d$-sphere $S^{\hspace{.2mm}\prime}$ by a bistellar move $\alpha
\mapsto \beta$ of index $\leq k-1$. Let $x$ be a vertex of $S$. If $x \not\in \alpha \sqcup\beta$ then ${\rm
lk}_S(x) = {\rm lk}_{S^{\hspace{.2mm}\prime}}(x)$ is $k$-stellated by induction hypothesis. If $x \in \alpha$
then ${\rm lk} _S(x)$ is obtained from the $k$-stellated sphere ${\rm lk}_{S^{\hspace{.2mm} \prime}}(x)$ by the
bistellar move $\alpha\setminus\{x\} \mapsto \beta$ of index $\leq k-1$. If $x \in \beta$ and $\beta\neq \{x\}$
then ${\rm lk} _S(x)$ is obtained from the $k$-stellated sphere ${\rm lk}_{S^{\hspace{.2mm}\prime}}(x)$ by the
bistellar move $\alpha \mapsto \beta\setminus\{x\}$ of index $\leq k-2$. If $\beta = \{x\}$ then ${\rm lk}_S(x)$
is the standard sphere $\partial \alpha$. Thus, in all cases, ${\rm lk}_S(x)$ is $k$-stellated. This proves part
$(a)$.

Let $S$ be a $k$-stacked $d$-sphere. Let $B$ be a $k$-stacked $(d+1)$-ball such that $\partial B = S$. If $x$ is
a vertex of $S$ then $x$ is a vertex of $B$ and $B^{\hspace{.2mm}\prime} = {\rm lk}_B(x)$ is a $d$-ball with
$\partial B^{\hspace{.2mm}\prime} = {\rm lk}_S(x)$. Therefore, it suffices to show that $B^{\hspace{.2mm}\prime}$
is also $k$-stacked. Indeed, if $\gamma$ is a face of $B^{\hspace{.2mm}\prime}$ of codimension $\geq k+1$ then
$\gamma \cup\{x\}$ is a face of $B$ of codimension $\geq k+1$, and hence $\gamma \cup \{x\} \in \partial B = S$,
so that $\gamma \in {\rm lk}_S(x) = \partial B^{\hspace{.2mm}\prime}$. \hfill $\Box$

\bigskip

The following is a stronger version of a result from \cite{bd13v2}. Also compare with Theorem 4.6 of \cite{mn2}.

\begin{theo} \label{theo-9}
Let $d \geq 2k+2$ and $M \in {\mathcal K}_k(d)$. Let $V(M)$ be the vertex set of $M$. Then \vspace{-3mm}
\begin{eqnarray} \label{eq2}
\overline{M} := \left\{\alpha \subseteq V(M) \, : \, \binom{\alpha}{\leq \, k+2\,} \subseteq M\right\}
\end{eqnarray}
is the unique homology $\,(d+ 1)$-manifold such that $\,\partial \overline{M} = M$ and $\,{\rm skel}_{d-
k}(\overline{M}) = {\rm skel}_{d- k}(M)$.
\end{theo}

\noindent {\bf Proof.} Fix $x \in V(\overline{M}) = V(M)$.

\smallskip

\noindent {\bf Claim\,:} ${\rm lk}_{\overline{M}}(x) = \overline{{\rm lk}_M(x})$, where the right hand side is as
defined in Theorem \ref{theo-6}.

\smallskip

From the definition, we see that $\alpha \in {\rm lk}_{\overline{M}}(x) \Rightarrow \alpha\sqcup\{x\}\in
\overline{M} \Rightarrow \binom{\alpha\sqcup \{x\}}{\leq k+2} \subseteq M \Rightarrow \binom{\alpha}{\leq k+1}
\subseteq {\rm lk}_M(x) \Rightarrow \alpha \in \overline{{\rm lk}_M(x})$. Thus, we have ${\rm
lk}_{\overline{M}}(x)\subseteq \overline{{\rm lk}_M(x})$.

Conversely, let $\alpha \in \overline{{\rm lk}_M(x})$. Then $\binom{\alpha}{\leq k+1} \subseteq {\rm lk}_M(x)$,
so that each $\gamma \subseteq \alpha \sqcup\{x\}$ such that $x \in \gamma$ and $\#(\gamma) \leq k+ 2$ is in $M$.
Therefore, to prove that $\alpha \in {\rm lk}_{\overline{M}}(x)$, it suffices to show that each $\gamma \subseteq
\alpha$ with $\#(\gamma) \leq k+2$ is in $M$. Since $\alpha \in \overline{{\rm lk}_M(x})$, such a set $\gamma$ is
in $\overline{{\rm lk}_M(x})$, and hence $\gamma \in {\rm skel}_{k+1}(\overline{{\rm lk}_{M}(x})) \subseteq {\rm
skel}_{d-k- 1}(\overline{{\rm lk}_{M}(x})) = {\rm skel}_{d-k- 1}({\rm lk}_{M}(x)) \subseteq {\rm lk}_M(x)
\subseteq M$. (Here the first inclusion holds since $k+1 \leq d-k-1$.) This proves that $\alpha \in
\overline{{\rm lk}_M(x}) \Rightarrow \alpha \in {\rm lk}_{\overline{M}}(x)$, so that $\overline{{\rm lk}_M(x})
\subseteq {\rm lk}_{\overline{M}}(x)$. This proves the claim.

\smallskip

In view of Theorem \ref{theo-6}, the claim implies that $\overline{M}$ is a homology $(d+1)$-manifold with
boundary, and ${\rm lk}_{\partial \overline{M}}(x) = \partial({\rm lk}_{\overline{M}}(x)) =
\partial(\overline{{\rm lk}_M(x})) = {\rm lk}_M(x)$ for every vertex $x$. Therefore, $\partial \overline{M} = M$,
and we have\,:
\begin{eqnarray*}
{\rm lk}_{{\rm skel}_{d-k}(\overline{M})}(x) & = & {\rm skel}_{d- k-1}({\rm lk}_{\overline{M}}(x)) = {\rm
skel}_{d-k- 1}(\overline{{\rm lk}_{M}(x})) = {\rm skel}_{d-k-1}({\rm lk}_{M}(x)) \\ &=&  {\rm lk}_{{\rm
skel}_{d-k}(M)}(x)
\end{eqnarray*}
for every vertex $x$. Thus, ${\rm skel}_{d-k}(\overline{M}) = {\rm skel}_{d-k}(M)$.

Now, if $N$ is any homology $(d+1)$-manifold with $\partial N = M$ and ${\rm skel}_{d-k}(N) = {\rm
skel}_{d-k}(M)$, then for any vertex $x$, we have\,:
\begin{eqnarray*}
\partial({\rm lk}_N(x)) & = & {\rm lk}_{\partial N}(x) = {\rm lk}_M(x), ~ \mbox{ and} \\
{\rm skel}_{d-k-1}({\rm lk}_N(x)) & = & {\rm lk}_{{\rm skel}_{d-k}(N)}(x) = {\rm lk}_{{\rm skel}_{d-k}(M)}(x)  =
{\rm skel}_{d-k-1}({\rm lk}_{M}(x)).
\end{eqnarray*}
Therefore, the uniqueness assertion in Theorem \ref{theo-6} implies that ${\rm lk}_N(x) = \overline{{\rm
lk}_M(x}) = {\rm lk}_{\overline{M}}(x)$ for every vertex $x$ and hence $N = \overline{M}$. This completes the
proof. \hfill $\Box$

\begin{remark} \label{remark-1}
{\rm If $M$ is a $k$-stacked homology sphere of dimension $d \geq 2k+2$ then $M \in {\mathcal K}_k(d)$ by Theorem
\ref{theo-8}. In this case, the uniqueness statements in Theorems \ref{theo-6} and \ref{theo-9} show that the two
definitions of $\overline{M}$ (given in (1) and (2)) agree. Also, if we define $\overline{\mathcal K}_k(d+1)$ to
be the class of all $(d+1)$-dimensional simplicial complexes all whose vertex links are $k$-stacked homology
$d$-balls, then by Theorems \ref{theo-6} and \ref{theo-9}, for $d\geq 2k+2$, $M \mapsto \overline{M}$ is a
bijection from ${\mathcal K}_k(d)$ onto $\overline{\mathcal K}_k(d+1)$. The boundary map provides its inverse.}
\end{remark}

\begin{remark} \label{remark-2}
{\rm In view of Theorem \ref{theo-9} above and Theorem 4.4 of \cite{mn2}, for $M \in {\mathcal K}_k(d)$ with
$d\geq 2k+2$, we have $H_i(M;\mathbb{Z}) = \{0\}$ for $k+1 \leq i \leq d-k-1$ and $H_k(M;\mathbb{Z})$ is torsion
free. A special case of this result (with the extra assumption of 2-neighbourliness) was proved in \cite[Theorem
3.7 (d)]{bd16}. }
\end{remark}

Following \cite{mn2}, we introduce an extension of Definition \ref{k-stacked-ball}.

\begin{defn} \label{k-stacked-manifold}
{\rm A homology manifold $N$ with boundary is said to be {\em $k$-stacked} if all its faces of codimension $k+1$
are in its boundary. A homology $d$-manifold $M$ without boundary is said to be {\em $k$-stacked}, if there is a
$k$-stacked homology $(d+1)$-manifold $N$ (with boundary) such that $M = \partial N$. Thus, Theorem \ref{theo-9}
says that, for dimension $d\geq 2k+2$, all members of ${\mathcal K}_k(d)$ are $k$-stacked. }
\end{defn}

Recall from \cite{bd16} that a simplicial complex $X$ is said to be {\em tight} with respect to a given field
$\mathbb{F}$ (or {\em $\mathbb{F}$-tight}) if, for every induced subcomplex $Y$ of $X$, the inclusion map $Y
\hookrightarrow X$ induces an injective morphism $H_j(Y; \mathbb{F}) \hookrightarrow H_j(X; \mathbb{F})$ for each
$j\geq 0$. In \cite{bd16} we proved that, for $d\neq 2k+1$, any $(k+1)$-neighborly $\mathbb{F}$-orientable member
of ${\mathcal W}_k(d)$ is $\mathbb{F}$-tight. Here we prove the following combinatorial characterization of
$\mathbb{F}$-tightness which covers the case $d=2k+1$. (Recall from \cite{bd16} that any $\mathbb{F}$-tight
homology manifold without boundary is $\mathbb{F}$-orientable.)

\begin{theo} \label{theo-10}
Let $M \in {\mathcal W}_k(2k+1)$ be $\mathbb{F}$-orientable and $(k+1)$-neighborly. Then the following are
equivalent\vspace{-1mm}
\begin{enumerate}
\item[$(i)$] $M$ is $\mathbb{F}$-tight, \vspace{-1mm} \item[$(ii)$] $M$ is $k$-stacked, and \vspace{-1mm}
\item[$(iii)$] $\beta_k(M; \mathbb{F}\hspace{.5mm}) = {\binom{n-k-3}{k +1}}/{\binom{2k+3}{k+1}}$, where $n =
f_0(M)$.
\end{enumerate}
\end{theo}

\noindent {\bf Proof.} $(i)\Leftrightarrow (iii)$ was proved in \cite[Theorem 3.10 (b)]{bd16}.

$(ii)\Rightarrow (iii)$ is immediate from \cite[Theorem 3.1]{mn2} applied to any $k$-stacked homology
$(2k+2)$-manifold $\Delta$ such that $\partial \Delta = M$. (Notice that, since $M$ has $n$ vertices and is
$(k+1)$-neighborly, the same is true for $\Delta$. Therefore, the $h$-vector of $\Delta$ satisfies
$h_{k+1}(\Delta) = \binom{n-k-3}{k+1}$.) Therefore, to complete the proof, it suffices to show that
$(i)\Rightarrow (ii)$.

So, let $M$ be $\mathbb{F}$-tight. Define $\overline{M}$ by Equation (\ref{eq2}). It suffices to show that
$\overline{M}$ is a $k$-stacked homology $(2k+2)$-manifold with boundary and $\partial \overline{M} = M$. As in
the proof of Theorem \ref{theo-9}, this will follow once we show that ${\rm lk}_{\overline{M}}(x) =
\overline{{\rm lk}_M(x)}$ for every $x\in V(M) = V(\overline{M})$. (Notice that ${\rm lk}_M(x)$ is a
$k$-stellated sphere of dimension $2k$, and hence Theorems \ref{theo-4} and \ref{theo-6} apply to it.) As in the
proof of Theorem \ref{theo-9}, the proof of ${\rm lk}_{\overline{M}}(x) \subseteq \overline{{\rm lk}_M(x)}$ is
easy. Also, to prove the reverse inclusion, it suffices to show that, whenever $\binom{\alpha}{\leq k+1}
\subseteq {\rm lk}_M(x)$ and $\gamma \in \binom{\alpha}{\leq k+2}$, we have $\gamma\in M$. In proving this, we
may assume without loss of generality that $\#(\alpha) = k+2$ and $\gamma = \alpha$. So, we are reduced to
proving that every missing $(k+1)$-face $\alpha$ of ${\rm lk}_M(x)$ belongs to $M$. Suppose not. Then the induced
subcomplex of $M$ on the vertex set $\alpha\sqcup\{x\}$ is the triangulated $(k+1)$-ball $B^0_1\ast
S^{\hspace{.15mm}k}_{k +2}$. Clearly, every induced subcomplex of an $\mathbb{F}$-tight simplicial complex is
$\mathbb{F}$-tight. Since $M$ is assumed to be $\mathbb{F}$-tight, it follows that the cone $B^0_1\ast
S^{\hspace{.15mm}k}_{k +2}$ over the standard sphere $S^{\hspace{.15mm}k}_{k +2}$ is $\mathbb{F}$-tight. This is
a contradiction since it is easy to see that the standard ball $B^{\hspace{.15mm}d}_{d +1}$ is the only
$\mathbb{F}$-tight $d$-ball. \hfill $\Box$

\section{Examples, counterexamples and questions}

\begin{eg}[Stellated versus stacked spheres]
\label{example-1}
\begin{itemize}
\item[{\bf (a)}] {\rm Let $S^{\hspace{.2mm}d}_{2d+2} = (S^{\hspace{.2mm}0}_{2})^{\ast\,(d+1)}$, the join of $d+1$
copies of $S^{\hspace{.2mm}0}_{2}$. Being the boundary complex of the $(d+1)$-dimensional cross polytope,
$S^{\hspace{.2mm}d}_{2d+2}$ is a polytopal $d$-sphere. By Theorem \ref{theo-5}, all polytopal $d$-spheres are
$d$-stellated. Therefore, $S^{\hspace{.2mm}d}_{2d+2}$ is $d$-stellated. Also, by Theorem \ref{theo-1}, it is
$d$-stacked. Since $S^{\hspace{.2mm}d}_{2d+2}$ is the clique complex of its edge graph (1-skeleton), it is not
$(d- 1)$-stacked. (If there was a $(d-1)$-stacked $(d+1)$-ball $B$ such that $\partial B =
S^{\hspace{.2mm}d}_{2d+2}$, then all the faces of $B$ would be cliques of the edge graph of
$S^{\hspace{.2mm}d}_{2d +2}$. But, all such cliques are in $S^{\hspace{.2mm}d}_{2d+2}$ itself.) For the same
reason, $S^{\hspace{.2mm}d}_{2d+2}$ does not contain any induced standard sphere except $S^{\hspace{.2mm}0}_{2}$.
Therefore, it does not admit any bistellar move of index $\geq 2$. Hence $S^{\hspace{.2mm}d}_{2d+2}$ is not
$(d-1)$-stellated. (By the comment following Corollary \ref{coro-2}, any $k$-stellated $d$-sphere, excepting
$S^{\hspace{.2mm}d}_{d +2}$, admits a bistellar move of index $>d-k$.)

}

\item[{\bf (b)}] {\rm It is more difficult to find examples of $(d+1)$-stellated $d$-spheres (i.e., combinatorial
$d$-spheres) which are not $d$-stellated. The following example is due to Dougherty, Faber and Murphy \cite{dfm}.

Let $S^{\,3}_{16}$ be the pure 3-dimensional simplicial complex with vertex set $\mathbb{Z}_{16} =
\mathbb{Z}/16\mathbb{Z}$ and an automorphism $i \mapsto i+1$ (mod 16). Modulo this automorphism, the basic facets
of $S^{\,3}_{16}$ are\,:
$$
\{0,1,4,6\}, \{0,1,4,9\}, \{0,1,6,14\}, \{0,1,8,9\}, \{0,1,8,10\}, \{0,1,10,14\}, \{0,2,9,13\}.
$$
Of these, the fourth facet generates an orbit of length 8, while each of the other facets generates an orbit of
length 16. Thus, $S^{\,3}_{16}$ has $1 \times 8 + 6 \times 16 = 104$ facets. The face vector of $S^{\,3}_{16}$ is
$(16, 120, 208, 104)$. Since $120 = \binom{16}{2}$, $S^{\,3}_{16}$ is 2-neighbourly and hence it does not allow
any bistellar 1-move. Also, it is easy to verify that $S^{\,3}_{16}$ has no edge of (minimum) degree 3 (and hence
it has no vertex of degree 4), so that it does not allow any bistellar move of index 2 or 3 either. (So,
$S^{\,3}_{16}$ is an {\em unflippable} $3$-sphere in the sense of \cite{dfm}\,: it does not allow any bistellar
move of positive index.) Thus, $S^{\,3}_{16}$ is not 3-stellated. (Being a combinatorial 3-sphere, it is of
course 4-stellated.) Following the proof of Theorem \ref{theo-1}, fix a vertex $x$ of $S^{\,3}_{16}$, and let
$B^{\,4}_{16} = \{\{x\} \sqcup\alpha \, : \, x \not\in \alpha\in S^{\,3}_{16}\}$. Then $ B^{\,4}_{16}$ is a
4-ball with $\partial B^{\,4}_{16} = S^{\,3}_{16}$. Since $S^{\,3}_{16}$ is 2-neighbourly, $ B^{\,4}_{16}$ is a
2-stacked ball, and hence $S^{\,3}_{16}$ is an example of a 2-stacked 3-sphere which is not even 3-stellated. If
$ B^{\,4}_{16}$ was shellable, then (by Theorem \ref{theo-2}) it would be 2-shelled and hence (by Corollary
\ref{coro-1}) $S^{\,3}_{16}$ would be 2-stellated. Thus, $ B^{\,4}_{16}$ is an example of a non-shellable
2-stacked 4-ball. }

\item[{\bf (c)}] {\rm It is even more difficult to find examples of homology $d$-spheres which are not
$(d+1)$-stellated (i.e., not combinatorial $d$-spheres). Trivially, all homology spheres of dimension $d\leq 2$
are combinatorial spheres. In consequence, all homology manifolds of dimension $d\leq 3$ are combinatorial
manifolds. In \cite{ed} and \cite{f}, Edwards and Freedman proved that a triangulated homology manifold of
dimension $d\geq 3$ is a triangulated manifold if and only if all its vertex links are simply connected. In
conjunction with Perelman's theorem (3-dimensional Poincar\'{e} conjecture) this shows that all triangulated
4-manifolds are combinatorial manifolds. The (non-) existence of triangulated 4-spheres which are not
combinatorial spheres is equivalent to the still unresolved 4-dimensional smooth Poincar\'{e} conjecture.
(According to \cite{bd5}, any such 4-sphere would require at least 13 vertices.) Thus, $d=5$ is the smallest
dimension in which we may reasonably expect triangulated spheres which are not combinatorial spheres. The
following 16-vertex triangulation $\Sigma^{\,3}_{16}$ of the Poincar\'{e} (integral) homology 3-sphere was found
by Bj\"{o}rner and Lutz \cite{bl}. The vertices of $\Sigma^{\,3}_{16}$ are $1, \dots, 9, 1^{\hspace{.1mm}\prime},
\dots, 7^{\hspace{.2mm}\prime}$. Its facets are\,: $1  2  4  9$, $1  2  4  6^{\hspace{.2mm}\prime} $,
$1265^{\hspace{.2mm}\prime} $, $1266^{\hspace{.2mm}\prime}$, $1295^{\hspace{.2mm}\prime} $,
$1343^{\hspace{.2mm}\prime}$, $1346^{\hspace{.2mm}\prime}$, $1371^{\hspace{.1mm}\prime} $,
$1373^{\hspace{.2mm}\prime} $, $ 1 3 1^{\hspace{.1mm}\prime} 6^{\hspace{.2mm}\prime} $, $1  4 9
3^{\hspace{.2mm}\prime} $, $1564^{\hspace{.2mm}\prime} $, $ 15  6 5^{\hspace{.2mm}\prime} $, $158
2^{\hspace{.2mm}\prime} $, $ 1584^{\hspace{.2mm}\prime} $, $152^{\hspace{.2mm}\prime}5^{\hspace{.2mm}\prime}$,
$164^{\hspace{.2mm}\prime}6^{\hspace{.2mm}\prime}$, $1781^{\hspace{.1mm}\prime}$, $1782^{\hspace{.2mm}\prime}$,
$172^{\hspace{.2mm}\prime}3^{\hspace{.2mm}\prime}$, $181^{\hspace{.1mm}\prime}4^{\hspace{.2mm}\prime}$,
$192^{\hspace{.2mm}\prime}3^{\hspace{.2mm}\prime}$, $192^{\hspace{.2mm}\prime}5^{\hspace{.2mm}\prime}$,
$11^{\hspace{.1mm}\prime}4^{\hspace{.2mm}\prime}6^{\hspace{.2mm}\prime}$, $2351^{\hspace{.1mm}\prime}$,
$2352^{\hspace{.2mm}\prime}$, $2371^{\hspace{.1mm}\prime}$, $2374^{\hspace{.2mm}\prime}$,
$232^{\hspace{.2mm}\prime}4^{\hspace{.2mm}\prime}$, $2494^{\hspace{.2mm}\prime}$,
$242^{\hspace{.2mm}\prime}4^{\hspace{.2mm}\prime}$, $242^{\hspace{.2mm}\prime}6^{\hspace{.2mm}\prime}$,
$2582^{\hspace{.2mm}\prime} $, $ 2 5 8  3^{\hspace{.2mm}\prime} $, $ 2 5 1^{\hspace{.1mm}\prime}
3^{\hspace{.2mm}\prime} $, $261^{\hspace{.1mm}\prime} 3^{\hspace{.2mm}\prime} $, $261^{\hspace{.1mm}\prime}
5^{\hspace{.2mm}\prime} $, $263^{\hspace{.2mm}\prime} 6^{\hspace{.2mm}\prime} $, $2794^{\hspace{.2mm}\prime}$, $
2 7  9  5^{\hspace{.2mm}\prime}$, $ 2 7 1^{\hspace{.1mm}\prime} 5^{\hspace{.2mm}\prime} $,
$282^{\hspace{.2mm}\prime} 6^{\hspace{.2mm}\prime} $, $283^{\hspace{.2mm}\prime} 6^{\hspace{.2mm}\prime} $,
$3455^{\hspace{.2mm}\prime}$, $ 3 4  5  6^{\hspace{.2mm}\prime} $, $ 3  4 3^{\hspace{.2mm}\prime}
5^{\hspace{.2mm}\prime} $, $351^{\hspace{.1mm}\prime} 6^{\hspace{.2mm}\prime} $, $352^{\hspace{.2mm}\prime}
5^{\hspace{.2mm}\prime} $, $373^{\hspace{.2mm}\prime} 4^{\hspace{.2mm}\prime} $,
$32^{\hspace{.2mm}\prime}4^{\hspace{.2mm}\prime}5^{\hspace{.2mm}\prime}$,
$33^{\hspace{.2mm}\prime}4^{\hspace{.2mm}\prime}5^{\hspace{.2mm}\prime}$, $    4  5  6  7 $, $ 4  5  6
5^{\hspace{.2mm}\prime} $, $4576^{\hspace{.2mm}\prime}$, $4 6  7  2^{\hspace{.2mm}\prime} $, $ 4 6
1^{\hspace{.1mm}\prime} 2^{\hspace{.2mm}\prime} $, $461^{\hspace{.1mm}\prime} 5^{\hspace{.2mm}\prime} $,
$472^{\hspace{.2mm}\prime} 6^{\hspace{.2mm}\prime} $, $4893^{\hspace{.2mm}\prime}$, $4894^{\hspace{.2mm}\prime}$,
$481^{\hspace{.1mm}\prime} 4^{\hspace{.2mm}\prime} $, $481^{\hspace{.1mm}\prime} 5^{\hspace{.2mm}\prime} $,
$483^{\hspace{.2mm}\prime} 5^{\hspace{.2mm}\prime} $,
$41^{\hspace{.1mm}\prime}2^{\hspace{.2mm}\prime}4^{\hspace{.2mm}\prime}$, $5674^{\hspace{.2mm}\prime} $,
$5794^{\hspace{.2mm}\prime}$, $5796^{\hspace{.2mm}\prime} $, $5893^{\hspace{.2mm}\prime}$,
$5894^{\hspace{.2mm}\prime}$, $591^{\hspace{.1mm}\prime}3^{\hspace{.2mm}\prime}$, $591^{\hspace{.1mm}\prime}
6^{\hspace{.2mm}\prime}$, $672^{\hspace{.2mm}\prime} 3^{\hspace{.2mm}\prime}$, $673^{\hspace{.2mm}\prime}
4^{\hspace{.2mm}\prime}$, $61^{\hspace{.1mm}\prime}2^{\hspace{.2mm}\prime}3^{\hspace{.2mm}\prime}$,
$63^{\hspace{.2mm}\prime}4^{\hspace{.2mm}\prime}6^{\hspace{.2mm}\prime}$,
$781^{\hspace{.1mm}\prime}5^{\hspace{.2mm}\prime}$, $782^{\hspace{.2mm}\prime} 6^{\hspace{.2mm}\prime}$,
$785^{\hspace{.2mm}\prime} 6^{\hspace{.2mm}\prime}$, $795^{\hspace{.2mm}\prime} 6^{\hspace{.2mm}\prime}$,
$83^{\hspace{.2mm}\prime}5^{\hspace{.2mm}\prime}6^{\hspace{.2mm}\prime}$,
$91^{\hspace{.1mm}\prime}2^{\hspace{.2mm}\prime}3^{\hspace{.2mm}\prime}$,
$91^{\hspace{.1mm}\prime}2^{\hspace{.2mm}\prime}7^{\hspace{.2mm}\prime}$,
$91^{\hspace{.1mm}\prime}6^{\hspace{.2mm}\prime}7^{\hspace{.2mm}\prime}$,
$92^{\hspace{.2mm}\prime}5^{\hspace{.2mm}\prime}7^{\hspace{.2mm}\prime}$,
$95^{\hspace{.2mm}\prime}6^{\hspace{.2mm}\prime}7^{\hspace{.2mm}\prime}$,
$1^{\hspace{.1mm}\prime}2^{\hspace{.2mm}\prime} 4^{\hspace{.2mm}\prime}7^{\hspace{.2mm}\prime}$,
$1^{\hspace{.1mm}\prime}4^{\hspace{.2mm}\prime} 6^{\hspace{.2mm}\prime}7^{\hspace{.2mm}\prime}$,
$2^{\hspace{.2mm}\prime}4^{\hspace{.2mm}\prime} 5^{\hspace{.2mm}\prime}7^{\hspace{.2mm}\prime} $,
$3^{\hspace{.2mm}\prime}4^{\hspace{.2mm}\prime} 5^{\hspace{.2mm}\prime}6^{\hspace{.2mm}\prime} $,
$4^{\hspace{.2mm}\prime}5^{\hspace{.2mm}\prime} 6^{\hspace{.2mm}\prime}7^{\hspace{.2mm}\prime}$. The face vector
of $\Sigma^{\,3}_{16}$ is $(16,106, 180, 90)$. Bj\"{o}rner and Lutz conjectured that it is strongly minimal in
the sense that it has the componentwise minimum face vector among all possible triangulations of the Poincar\'e
homology sphere.

Note that the vertex $6^{\hspace{.2mm}\prime}$ is adjacent with all other vertices in $\Sigma^{\,3}_{16}$. Let
$D^{\,4}_{16}$ be the 4-dimensional simplicial complex whose facets are $\alpha \cup
\{6^{\hspace{.2mm}\prime}\}$, where $\alpha$ ranges over all facets of $\Sigma^{\,3}_{16}$ not containing the
vertex $6^{\hspace{.2mm} \prime}$. Define $S^{\,5}_{18}= \partial(D^{\,4}_{16} \ast B^{1}_{2})$, the boundary of
the join of $D^{\,4}_{16}$ and an edge. Observe that, $|S^{\,5}_{18}|$ is the double suspension of the
Poincar\'{e} homology sphere $|\Sigma^{\,3}_{16}|$. Therefore, by Cannon's double suspension theorem (cf.
\cite{ca}, actually Cannon's theorem is a straightforward consequence of the result of Edwards and Freedman
quoted above), $S^{\,5}_{18}$ is a triangulated 5-sphere. Since it has $\Sigma^{\,3}_{16}$ as the link of an
edge, $S^{\,5}_{18}$ is not a combinatorial sphere.

Let $D^{\,6}_{18} = D^{\,4}_{16} \ast B^{1}_{2}$, $D^{\,7}_{19} = D^{\,4}_{16} \ast B^{\hspace{.1mm}2}_{3}$ and
$S^{\,6}_{19} = \partial D^{\,7}_{19}$. By the above logic, $S^{\,6}_{19}$ is a triangulated 6-sphere. Since
$D^{\,6}_{18}$ is the antistar of a vertex in $S^{\,6}_{19}$, it follows from Lemma 4.1 in \cite{bd9} that
$D^{\,6}_{18}$ is a triangulated $6$-ball. Since the vertex $6^{\hspace{.2mm}\prime}$ is adjacent to all the
vertices in $\Sigma^{\,3}_{16}$, the construction of $D^{\,6}_{18}$ shows that all the 3-faces of $D^{\,6}_{18}$
lie in its boundary. Thus, $D^{\,6}_{18}$ is a 2-stacked triangulated ball. As $S^{\,5}_{18} = \partial
D^{\,6}_{18}$, it follows that $S^{\,5}_{18}$ is an example of 2-stacked 5-sphere which is not even 6-stellated.
(In \cite[Example 6.2]{mn}, Murai and Nevo present a series of 2-stacked balls which are not shellable.
Therefore, by Theorem \ref{theo-4}, their boundaries are examples of 2-stacked spheres which are not 2-stellated.
However, these examples seem to be 3-stellated.) }

\item[{\bf (d)}] {\rm Let $S$ be a triangulated $d$-sphere and $B$ be a $k$-stacked ball such that $\partial B =
S$. Then, for any $e\geq 0$, $B \ast B^{\hspace{.2mm}e}_{e+1}$ is a $k$-stacked ball, and hence $\partial(B \ast
B^{\hspace{.2mm}e}_{e+1})$ is a $k$-stacked $(d+e+1)$-sphere. Also, $S$ is a combinatorial sphere if and only if
$\partial(B \ast B^{\hspace{.2mm}e}_{e +1})$ is so. Applying this construction to the pair
$(S^{\hspace{.2mm}5}_{18}, D^{\hspace{.2mm}6}_{18})$ in example (c) above, we find that for each $d \geq 5$,
there are 2-stacked triangulated $d$-spheres which are not even $(d+1)$-stellated.

\noindent {\bf Claim.} If $B^{\hspace{.2mm}4}_{16}$ is as in example (b) above then
$\partial(B^{\hspace{.2mm}4}_{16} \ast B^{\hspace{.1mm}e}_{e+1})$ is unflippable.

For $e \geq 0$, let $\widetilde{B}^{\,e+5}_{e+17} : = B^{\,4}_{16} \ast B^{\,e}_{e+1}$ and $\widetilde{S}^{\,e
+4}_{e+17} := \partial \widetilde{B}^{\,e+5}_{e+17}$. Thus, $\widetilde{S}^{\,e+4}_{e+17} = (S^{\,3}_{16} \ast
B^{\,e}_{e+1}) \cup (B^{\,4}_{16} \ast S^{\,e-1}_{e+1})$. Since $S^{\,3}_{16}$ is 2-neighbourly, so is
$\widetilde{S}^{\,e +4}_{e+17}$. Therefore, $\widetilde{S}^{\,e +4}_{e+17}$ does not admit any bistellar 1-move.
Suppose, if possible, that $\alpha \mapsto \beta$ is a bistellar move of index $\geq 2$ on $\widetilde{S}^{\,e
+4}_{e + 17}$. Thus, ${\rm lk}_{\widetilde{S}^{\,e +4}_{e+17}}(\alpha) = \partial \beta$ and $\dim(\beta) \geq
2$, $\beta \not\in \widetilde{S}^{\,e +4}_{e +17}$. Write $\alpha = \alpha_1\sqcup \alpha_2$, where $\alpha_1$ is
a face of $B^{\,4}_{16}$ and $\alpha_2$ is a face of $B^{\,e}_{e +1}$. If $\alpha_1$ is an interior face of
$B^{\,4}_{16}$, then $\alpha_2 \in S^{\,e- 1}_{e+1}$ and $\partial \beta = {\rm lk}_{\widetilde{S}^{\,e
+4}_{e+17}}(\alpha) =  {\rm lk}_{B^{\,4}_{16}}(\alpha_1)\ast {\rm lk}_{S^{\,e-1}_{e+ 1}}(\alpha_2)$. Since the
standard sphere $\partial \beta$ can't be written as the join of two spheres, it follows that either $\alpha_1$
is a facet of $B^{\,4}_{16}$ or $\alpha_2$ is a facet of $S^{\,e-1}_{e+1}$. If $\alpha_1$ is a facet of
$B^{\,4}_{16}$, then $\partial \beta = {\rm lk}_{S^{\,e -1}_{e+1}}(\alpha_2)$ and hence $\beta \in B^{\,e}_{e+1}
\subseteq \widetilde{S}^{\,e +4}_{e + 17}$. This is a contradiction since $\overline{\alpha} \ast \partial \beta$
is an induced subcomplex of $\widetilde{S}^{\,e +4}_{e+17}$. So, $\alpha_2$ is a facet of $S^{\,e-1}_{e+1}$ and
hence $\partial \beta = {\rm lk}_{B^{\,4}_{16}}(\alpha_1)$. Since $B^{\hspace{.2mm}4}_{16}$ is a 2-stacked 4-ball
and $\alpha_1$ is an interior face of $B^{\hspace{.2mm}4}_{16}$, we get $2 \leq \dim(\alpha_1)= 4 - \dim(\beta)
\leq 2$ and hence $\dim(\alpha_1) = \dim(\beta) =2$. Let $\alpha_1= xuv$ (where $x$ is the fixed vertex chosen in
$S^{\,3}_{16}$ to construct $B^{\,4}_{16}$). Then ${\rm lk}_{S^{\,3}_{16}}(uv) = {\rm
lk}_{B^{\,4}_{16}}(\alpha_1) = \partial\beta$. This is not possible since $S^{\,3}_{16}$ does not contain any
edge of degree 3.

Thus $\alpha_1$ is a boundary face of $B^{\,4}_{16}$, i.e., $\alpha_1 \in S^{\,3}_{16}$. If $\alpha_2$ is the
facet of $B^{\,e }_{e+1}$ then ${\rm lk}_{S^{\,3}_{16}}(\alpha_1) = {\rm lk}_{\widetilde{S}^{\,e +4}_{e +
17}}(\alpha) = \partial \beta$. Hence $\dim(\alpha_1) \geq 2$ and therefore $\dim(\beta) \leq 1$, a
contradiction. So, $\alpha_2$ is not the facet of $B^{\,e}_{e+1}$ (and hence ${\rm lk}_{B^{\,e}_{e+1}}(\alpha_2)$
is a standard ball). Thus, the ball $B_1 := {\rm lk}_{B^{\,4}_{16}}(\alpha_1) \ast {\rm lk}_{B^{\,e}_{e+
1}}(\alpha_2)$ is a non-trivial join of balls, so that all the vertices of $B_1$ are in its boundary. But,
$\partial B_1 = {\rm lk}_{\widetilde{S}^{\,e+4}_{e+17}}(\alpha) = \partial \beta$. Therefore, $B_1$ is the
standard ball $\overline{\beta}$ and hence ${\rm lk}_{B^{\,4}_{16}}(\alpha_1)$ is a standard ball. Therefore,
${\rm lk}_{S^{\,3}_{16}}(\alpha_1)$ is a standard sphere and hence $\dim(\alpha_1) \geq 2$. So, ${\rm
lk}_{B^{\,4}_{16}}(\alpha_1)$ is a standard ball of dimension $\leq 1$, i.e., it is a vertex or an edge. Then the
vertex set of ${\rm lk}_{B^{\,4}_{16}}(\alpha_1)$ is a face in $S^{\,3}_{16}$. So, the vertex set $\beta$ of
${\rm lk}_{\widetilde{S}^{\,e +4}_{e+17}}(\alpha)$ is a face of $S^{\,3}_{16} \ast B^{\,e}_{e+1} \subseteq
\widetilde{S}^{\,e+4}_{e+17}$. Therefore, $\overline{\alpha} \ast \partial\beta$ is not an induced subcomplex of
$\widetilde{S}^{\,e+4}_{e+17}$, a contradiction. Thus, for each $e\geq 0$, $\widetilde{S}^{\,e+4}_{e+17}$ is an
unflippable combinatorial $(e+4)$-sphere.

From this claim, it follows that $\partial( B^{\,4}_{16} \ast B^{\hspace{.2mm}d-4}_{d-3})$ is a combinatorial
$d$-sphere which is not $d$-stellated. Since $B^{\,4}_{16}$ is a 2-stacked ball, it follows that $B^{\,4}_{16}
\ast B^{\hspace{.2mm}d-4}_{d -3}$ is also 2-stacked. This implies that $\partial( B^{\,4}_{16} \ast
B^{\hspace{.2mm}d-4}_{d-3})$ is a 2-stacked combinatorial $d$-sphere which is not $d$-stellated, for $d\geq 4$.
From this and the observation in (b), we find that for each $d\geq 3$, there are $2$-stacked combinatorial
$d$-spheres which are not $d$-stellated.

Since the classes $\Sigma_k(d)$, ${\mathcal S}_k(d)$ are increasing in $k$, we get\,:}
\begin{itemize}
\item[$\bullet$] For all $(k, l, d)$ such that $2 \leq k \leq l \leq d$ and $d\geq 3$,  there are $k$-stacked
combinatorial $d$-spheres which are not $l$-stellated.

\item[$\bullet$] For all $(k, l, d)$ such that $2 \leq k \leq l \leq d+1$ and $d\geq 5$,  there are $k$-stacked
triangulated $d$-spheres which are not $l$-stellated.
\end{itemize}
\item[{\bf (e)}] {\rm Let $S^{\hspace{.2mm}3}_{10}$ be the pure simplicial complex of dimension three whose
vertices are the digits $0, 1, \dots, 9$ and whose facets are\,:
\begin{eqnarray*}
0123, 1234, 2345, 3456, 4567, 5678, 6789, 0128, 0139, 0189, 0238, 0356, 0358, 0369,  \\
0568, 0689, 1248, 1349, 1457, 1458, 1467, 1469, 1578, 1679, 1789, 2358, 2458, 3469.
\end{eqnarray*}
Let $S^{\hspace{.2mm}2}_{10}$ be the pure 2-dimensional subcomplex of $S^{\hspace{.2mm}3}_{10}$ whose facets
are\,:
$$
012, 013, 023, 124, 134, 235, 245, 346, 356, 457, 467, 568, 578, 679, 689, 789.
$$
Clearly, $S^{\hspace{.2mm}2}_{10}$ is a triangulated 2-sphere.

Let $B_1$ (respectively, $B_2$) be the 3-dimensional subcomplex of $S^{\hspace{.2mm}3}_{10}$ whose facets are the
first seven (respectively, last twenty one) facets of $S^{\hspace{.2mm}3}_{10}$. Then $B_1$ is a normal
pseudomanifold and the dual graph of $B_1$ is a path. So, by Theorem \ref{theo-7}, $B_1$ is a 1-stacked 3-ball.
It is easy to see that $B_2$ is a triangulated 3-manifold with $\partial B_2 = S^{\hspace{.2mm}2}_{10} = \partial
B_1$. It is not difficult to check that $B_2$ is collapsible and hence a triangulated 3-ball. This implies that
$S^{\hspace{.2mm}3}_{10}$ is a triangulated 3-sphere. (So, $S^{\hspace{.2mm}2}_{10}$ is a triangulated 2-sphere
embedded in $S^{\hspace{.2mm}3}_{10}$ and divides $S^{\hspace{.2mm}3}_{10}$ into two closed ``hemispheres" $B_1$
and $B_2$.)

Since $B_1$ is a 1-stacked 3-ball and $S^{\hspace{.2mm}2}_{10} = \partial B_1$, $S^{\hspace{.2mm}2}_{10}$ is
1-stellated. But, $S^{\hspace{.2mm}2}_{10}$ also bounds the ball $B_2$ which is Ziegler's example \cite{z2} of a
non-shellable 3-ball\,! (If $\alpha$ is a facet of a triangulated $d$-ball $B$, then one says $\alpha$ is an {\em
ear} of $B$ if $B \setminus \{\alpha\}$ is also a triangulated $d$-ball. Clearly, if $B$ is shellable, then the
last facet, added while obtaining $B$ from $B^{\hspace{.2mm}d}_{d+1}$ by a sequence of shelling moves, must be an
ear of $B$. Thus, if $B$ has no ears, then it must be non-shellable. Such balls are ``strongly non-shellable" in
the terminology of Ziegler. A facet $\alpha$ of $B$ is an ear of $B$ if and only if the induced subcomplex of
$\partial B$ on the vertex set $\alpha$ is a $(d-1)$-ball. Using this criterion, it is possible to verify that
$B_2$ has no ears\,: it is strongly non-shellable.) }

\item[{\bf (f)}] {\rm The following example of a shellable 3-ball with a unique ear is due to Frank Lutz
\cite{lu3}. Consider the pure 3-dimensional 2-neighbourly simplicial complex $S^{\hspace{.2mm}3}_{8}$ with
vertices $1, 2, \dots, 8$ and facets
\begin{eqnarray*}
1234, 2345, 3456, 4567, 5678, 1237, 1248, 1278, 1348, 1356, \\
1357, 1368, 1568, 1578, 2357, 2457, 2467, 2468, 2678, 3468.
\end{eqnarray*}
Let $S^{\hspace{.2mm}2}_{8}$ be the pure 2-dimensional subcomplex of $S^{\hspace{.2mm}3}_{8}$ with facets
$$
123, 124, 134, 235, 245, 346, 356, 457, 467, 568, 578, 678.
$$
Again, $S^{\hspace{.2mm}2}_{8}$ is a triangulated 2-sphere embedded in the triangulated 3-sphere
$S^{\hspace{.2mm}3}_{8}$. As in (e) above, $S^{\hspace{.2mm}2}_{8}$ divides $S^{\hspace{.2mm}3}_{8}$ into two
3-balls $B_1$ and $B_2$. The facets of $B_1$ are the first five facets of $S^{\hspace{.2mm}3}_{8}$, while the
facets of $B_2$ are the remaining fifteen facets of $S^{\hspace{.2mm}3}_{8}$. Again, $B_1$ is a 1-stacked 3-ball
since its dual graph is a path. We have $\partial B_1 = S^{\hspace{.2mm}2}_{8} = \partial B_2$. Thus,
$S^{\hspace{.2mm}2}_{8}$ is a 1-stellated sphere. The other ball $B_2$ bounded by $S^{\hspace{.2mm}2}_{8}$ is
shellable (indeed, 2-shelled). (A shelling of $B_2$\,: 1357, 1356, 1368, 1348, 1248, 3468, 1568, 1578, 1278,
2468, 2678, 1237, 2467, 2357, 2457.) But, $B_2$ has only one ear, namely $2457$.

Clearly, the class ${\mathcal S}_k(d)$ of $k$-stacked $d$-spheres is closed under connected sum. In consequence,
the class $\Sigma_1(d)$ of 1-stellated $d$-spheres is closed under connected sums. However, consider the
following construction. Take a standard 2-ball $B^{\hspace{.2mm}2}_{3}$ with a vertex set $\{a, b, c\}$ disjoint
from $V(S^{\hspace{.2mm}3}_{8})$, and form the join $B := B_2 \ast B^{\hspace{.2mm}2}_{3}$. Then $B$ is a
2-shelled 6-ball with a unique ear $2457abc$. Thus, $S := \partial B$ is a 2-stellated 5-sphere. The facets
$245abc$, $457abc$ are two of the facets of $S$ in the unique ear of $B$. Take a vertex disjoint copy
$B^{\hspace{.2mm}\prime}$ of $B$, and let $S^{\hspace{.2mm}\prime} = \partial B^{\hspace{.2mm}\prime}$, the
corresponding copy of $S$. Let $1^{\hspace{.1mm}\prime}, \dots, 8^{\hspace{.2mm}\prime}, a^{\hspace{.1mm}\prime},
b^{\hspace{.1mm}\prime}, c^{\hspace{.1mm}\prime}$ be the vertices of $B^{\hspace{.2mm}\prime}$ corresponding to
the vertices $1, \dots, 8, a, b, c$ respectively. Form the connected sum $\widetilde{B} = B \# B^{\hspace{.2mm}
\prime}$ by doing the identifications $2 \equiv 2^{\hspace{.2mm}\prime}, 4 \equiv 4^{\hspace{.2mm}\prime}, 5
\equiv 5^{\hspace{.2mm}\prime}, a \equiv a^{\hspace{.2mm}\prime}, b \equiv b^{\hspace{.2mm}\prime}, c \equiv
c^{\hspace{.2mm}\prime}$. Then $\widetilde{B}$ is a 16-vertex non-shellable 2-stacked 6-ball. Let $\widetilde{S}
= \partial \widetilde{B}$. Then $\widetilde{S}$ is a 16-vertex 2-stacked 5-sphere which is not 2-stellated (by
Theorems \ref{theo-6} and \ref{theo-4}). (It can be shown that $\widetilde{S}$ is 5-stellated.) But,
$\widetilde{S} = S \# S^{\hspace{.2mm}\prime}$, the connected sum of two 2-stellated 5-spheres. For $d\geq 5$, if
we take $B^{\,d- 3}_{d-2}$ in place of $B^{\,2}_3$ in the above construction then, by the same argument, we get a
$d$-sphere which is not 2-stellated and is the connected sum of two 2-stellated $d$-spheres. Thus
\begin{itemize}
\item[$\bullet$] For $d\geq 5$, the class $\Sigma_2(d)$ is not closed under connected sum.
\end{itemize}

}
\end{itemize}
\end{eg}

By Theorem \ref{theo-4}, all the $k$-stellated spheres of dimension $d \geq 2k-1$ are $k$-stacked. But, we are so
far unable to answer\,:

\begin{qn} \label{question-1}
{\rm Is there a $k$-stellated $d$-sphere which is not $k$-stacked\,?}
\end{qn}
Note that, by Theorems \ref{theo-1} and \ref{theo-4}, for an affirmative answer to Question \ref{question-1}, we
must have $k+1 \leq d \leq 2k-2$, and hence $k\geq 3$, $d\geq 4$.

Recall that a triangulated sphere is said to be {\em polytopal\,} if it is isomorphic to the boundary complex of
a simplicial polytope. We pose\,:

\begin{conj} \label{conj-1}
For $d\geq 2k$, a polytopal $d$-sphere is $k$-stellated if $($and only if\,$)$ it is $k$-stacked. Equivalently
$($in view of Theorems $\ref{theo-2}$ and $\ref{theo-4})$, if $S$ is a $k$-stacked polytopal sphere of dimension
$d \geq 2k$, then the $(d+1)$-ball $\overline{S}$ $($given by formula $(\ref{eq1}))$ is shellable.
\end{conj}

\begin{eg}
\label{example-2} {\rm Let $S = S^{\hspace{.2mm}k-1}_{k+1} \ast S^{\hspace{.2mm}k- 1}_{k+1}$, $B_1 =
S^{\hspace{.2mm}k-1}_{k+1} \ast B^{\hspace{.2mm}k}_{k +1}$ and $B_2 = B^{\hspace{.2mm}k}_{k+1} \ast
S^{\hspace{.2mm}k-1}_{k+1}$. Then $B_1$, $B_2$ are $k$-stacked polytopal $2k$-balls with $\partial B_1 = S =
\partial B_2$. Therefore, $S$ is a $(2k-1)$-dimensional $k$-neighbourly polytopal $k$-stacked sphere. Hence $S$ is
$k$-stellated by Theorem \ref{theo-5}. Thus, $S$ is an example of a $(2k- 1)$-dimensional $k$-stellated polytopal
sphere which bounds two distinct (though isomorphic) $k$-stacked balls. So, the bound $d\geq 2k$ in Theorem
\ref{theo-6} is sharp. }
\end{eg}

\begin{eg}[The Klee-Novik construction]
\label{example-3} {\rm For $d\geq 1$, let $S^{\hspace{.2mm}d+1}_{2d+4}$ be the join of $d+2$ copies of
$S^{\hspace{.2mm}0}_{2}$ with disjoint vertex sets $\{x_i, y_i\}$, $1\leq i \leq d+2$. Then $S^{\hspace{.2mm}d
+1}_{2d+4}$ is a triangulated sphere with missing edges $x_iy_i$, $1\leq i \leq d+2$ (cf. Example
\ref{example-1}\,(a)). Each of the $2^{d+2}$ facets of $S^{\hspace{.2mm}d +1}_{2d+4}$ may be encoded by a
sequence of $d+2$ signs as follows. If $\sigma$ is a facet, then for each index $i$ ($1\leq i \leq d+2$) $\sigma$
contains either $x_i$ or $y_i$, but not both. Put $\varepsilon_i = +$ if $x_i\in \sigma$ and $\varepsilon_i = -$
if $y_i\in \sigma$. Thus the sign sequence $(\varepsilon_1, \dots, \varepsilon_{d+2})$ encodes the facet
$\sigma$. For $0\leq k\leq d$, let $\overline{M}(k, d)$ be the pure $(d+1)$-dimensional subcomplex of
$S^{\hspace{.2mm}d +1}_{2d+4}$ whose facets are those facets $\sigma$ (of the latter complex) whose sign
sequences have at most $k$ sign changes. (A sign change in the sign sequence $(\varepsilon_1, \dots,
\varepsilon_{d+2})$ is an index $1\leq i \leq d+1$ such that $\varepsilon_{i+1} \neq \varepsilon_{i}$.) Then
$\overline{M}(k, d)$ is a pseudomanifold with boundary. Klee and Novik \cite{kn} proved that $M(k, d) := \partial
\overline{M}(k, d)$ is a triangulation of $S^{\hspace{.2mm}k} \times S^{\hspace{.2mm}d-k}$ for $0\leq k \leq d$.
(In their paper, Klee and Novik use the notation $B(k, d+2)$ for $\overline{M}(k, d)$.) The authors of \cite{kn}
observed that the permutations $D$, $E$ and $R$ are automorphisms of $\overline{M}(k, d)$ (and hence of $M(k,
d)$), where $D = {\displaystyle \prod_{j=1}^{d+2}(x_j, y_j)}$, $E = \hspace{-3mm}{\displaystyle \prod_{1\leq j<
(d+3)/2}\hspace{-3mm} (x_j, x_{d+3-j})(y_j, y_{d+3-j})}$ and $R = (x_1, \dots, x_{d+2})(y_1, \dots, y_{d+2})$
when $k$ is even, $R = (x_1, \dots, x_{d+2}, y_1, \dots, y_{d+2})$ when $k$ is odd. Clearly, these three
automorphisms generate a vertex-transitive automorphism group of $\overline{M}(k, d)$. Therefore, the links in
$\overline{M}(k, d)$ (or in $M(k, d)$) of all the vertices are isomorphic. The involution $A = {\displaystyle
\prod_{j ~ {\rm even}}(x_j, y_j)}$ is an isomorphism between $M(k, d)$ and $M(d-k, d)$. Therefore, in discussing
these constructions we may (and do) assume $d\geq 2k$. (However, $A$ is not an isomorphism between
$\overline{M}(k, d)$ and $\overline{M}(d-k, d)$. Indeed, $A$ maps $\overline{M}(k, d)$ to the ``complement" of
$\overline{M}(d-k, d)$ in $S^{\hspace{.2mm}d+1}_{2d+4}$.)

Let $I = \{1, 2, \dots, d+1\}$. Define the linear order $\prec$ on $\binom{I}{\leq \,k\,}$ by\,: $\alpha \prec
\beta$ if either $\#(\alpha) < \#(\beta)$ or else $\#(\alpha) = \#(\beta)$, $\alpha <_{\rm lex} \beta$, where
$<_{\rm lex}$ is the usual lexicographic order. Let $L$ be the link of the vertex $x_{d+2}$ in $\overline{M}(k,
d)$. Clearly, for each $\alpha \in \binom{I}{\leq \,k\,}$, there is a unique facet $\tau_{\alpha}$ of $L$ such
that $\alpha$ is precisely the set of sign-changes corresponding to the facet $\tau_{\alpha}\cup \{x_{d+2}\}$ of
$\overline{M}(k, d)$. We may transfer the linear order $\prec$ to the set of facets of $L$ via the bijection
$\alpha \mapsto \tau_{\alpha}$. Then, Klee and Novik show in \cite{kn} that $\prec$ is a shelling order for $L$.
Thus, $L$ is a shellable $d$-ball. What is more, if $\#(\alpha) = j \leq k$ then the facet $\tau_{\alpha}$ of $L$
is obtained (from the $d$-ball with facets $\tau_{\beta}$, $\beta\prec \alpha$) by a shelling move of index
$j-1$. In consequence, $L$ is a $k$-shelled $d$-ball. Since the automorphism group of $\overline{M}(k, d)$ is
vertex transitive, it follows that all vertex links of $\overline{M}(k, d)$ are $k$-shelled $d$-balls. Thus,
$\overline{M}(k, d)$ is a $(d+1)$-manifold with boundary. Also, since the boundary of a $k$-shelled ball is a
$k$-stellated sphere (Corollary \ref{coro-1}), it follows that $M(k, d) = \partial \overline{M}(k, d)$ has
$k$-stellated vertex links. Thus,
\begin{itemize}
\item[$\bullet$] $M(k, d)\in {\mathcal W}_k(d)$ for $d\geq 2k$.
\end{itemize}
Also note that, when $d\geq 2k+1$, the vertex links of $\overline{M}(k, d)$ are the unique (Theorem \ref{theo-6})
$k$-stacked balls bounded by the corresponding vertex links of $M(k, d)$. Therefore, $\overline{M}(k, d)$ is the
unique $(d+1)$-manifold $\overline{M}$ such that $\partial \overline{M} = M(k, d)$ and ${\rm skel}_{d-
k}(\overline{M}) = {\rm skel}_{d-k}(M(k, d))$. (In consequence, when $d\geq 2k+2$, $\overline{M}(k, d)$ may be
recovered from $M(k, d)$ via the formula (\ref{eq2}) above; cf. Theorem \ref{theo-9}.) Therefore, for $d \geq
2k+1$, every automorphism of $M(k, d)$ extends to an automorphism of $\overline{M}(k, d)$\,: they have the same
automorphism group. However, it is elementary to verify that the full automorphism group of $\overline{M}(k, d)$
is of order $4d+8$. (Since this group is transitive on the $2d+4$ vertices, it suffices to show that the full
stabilizer of the vertex $x_{d+2}$ is of order 2. This is easy.) Thus,
\begin{itemize}
\item[$\bullet$] When $d\geq 2k+1$, the full automorphism group of $M(k, d)$ is of order $4d+8$ (namely, the
group generated by $D$, $E$, $R$ above).
\end{itemize}
This leaves open the following tantalizing question. }
\end{eg}

\begin{qn} \label{question-2}
{\rm What is the full automorphism group of $M(k, 2k)$\,?}
\end{qn}

Notice that the involution $A$ defined above is also an automorphism of $M(k, 2k)$. However, since $A$ maps
$\overline{M}(k, 2k)$ to its complement in $S^{\hspace{.2mm}d +1}_{2d+4}$, $A$ is not an automorphism of
$\overline{M}(k, 2k)$. Therefore, $A \not\in H := \langle D, E, R\rangle$. The automorphism $A$ normalizes $H$,
so that the group $G := \langle D, E, R, A\rangle$ is of order $2\times \#(H) = 16(k+1)$. We suspect that $G$ is
the full automorphism group of $M(k, 2k)$. Is it\,?

\bigskip

\noindent {\bf Acknowledgement\,:} The authors thank Frank H. Lutz for drawing their attention to the reference
\cite{dfm} and for providing the unique ear 3-ball $B_2$ of Example \ref{example-1}\,(f) \cite{lu3}. The authors
also thank the anonymous referees for some useful comments. The second author was partially supported by grants
from UGC Centre for Advanced Study.

\end{document}